\documentclass[10pt,a4paper]{article}
\usepackage{amsmath, amsthm, amsfonts, amssymb, mathrsfs}
\usepackage[english]{babel}
\usepackage[all]{xy}
\usepackage{hyperref}
\usepackage{cite}
\usepackage[title]{appendix}
\usepackage[textwidth=12.9cm,textheight=21.5cm]{geometry} % A4 = 210mmx297mm
\usepackage{cleveref}
\theoremstyle{definition}
\newtheorem{teo}{Theorem}[subsection]
\newtheorem{teointro}{Theorem}
\newtheorem{prop}[teo]{Proposition}
\newtheorem{cor}[teo]{Corollary}
\newtheorem{lemma}[teo]{Lemma}
\newtheorem{rem}{Remark}[section]
\newtheorem*{remo}{Remark}
\newtheorem{deff}[teo]{Definition}
\newtheorem*{teobis}{Theorem 5.1.4'}

\def\ZZ{\mathbb{Z}}
\def\CC{\mathbb{C}}
\def\PP{\mathbb{P}}

\def\UU{\mathbb{U}}
\def\NN{\mathbb{N}}

\def\OO{\mathcal{O}}

\def\FF{\mathcal{F}}
\def\LL{\mathcal{L}}
\def\RR{\mathcal{R}}
\def\CCal{\mathcal{C}}

\def\address#1{\noindent {#1}}

\title{\bf Unfoldings and Deformations of Rational and Logarithmic Foliations}
\author{Ariel Molinuevo\thanks{The author was fully supported by CONICET, Argentina.}\address{}}
\date{}

\begin{document}

\maketitle

\centerline{\it Departamento de Matem\'atica, FCEN, Universidad de Buenos Aires, Argentina.}

\

\begin{abstract}
We study codimension one foliations in projective space $\PP^n$ over $\CC$ by looking at its first order perturbations: unfoldings and deformations. We give special attention to foliations of rational and logarithmic type.

For a differential form $\omega$ defining a codimension one foliation, we present a graded module $\UU(\omega)$, related to the first order unfoldings of $\omega$. If $\omega$ is a generic form of rational or logarithmic type, as a first application of the construction of $\UU(\omega)$, we classify the first order deformations that arise from first order unfoldings. Then, we count the number of isolated points in the singular set of $\omega$, in terms of a Hilbert polynomial associated to $\UU(\omega)$.

We review the notion of regularity of $\omega$ in terms of a long complex of graded modules that we also introduce in this work. We use this complex to prove that, for generic rational and logarithmic foliations, $\omega$ is regular if and only if every unfolding is trivial up to isomorphism.
\end{abstract}

\section{Introduction}

\hspace{9pt}An algebraic foliation of codimension one and degree $e-2$ in projective space $\PP^n$ over $\CC$, is given by a global section $\omega$ of 
the sheaf of twisted differential 1-forms $\Omega^1_{\PP^n}(e)$ that verifies the Frobenius integrability condition $\omega\wedge d\omega = 0$.
The space of such foliations forms a projective variety $\FF^1(\PP^n)(e)$. 
For arbitrary $n$ and $e$ not much is known about the irreducible components of $\FF^1(\PP^n)(e)$. 
We refer the reader to \cite{jou} and \cite{celn} where they classify the space of foliations of degree 1 and 2 respectively in $\PP^n$. 

\

There are two natural ways to make a first order perturbation of a codimension one foliation defined by $\omega$, namely first order deformations and unfoldings. The first order deformations, are given by a family of differential forms $\omega_\varepsilon$ parametrized by an infinitesimal parameter $\varepsilon$, such that $\omega_\varepsilon$ is integrable for every fixed $\varepsilon$ and at the origin coincides with $\omega$ up to scalar multiplication.
On the other side, first order unfoldings have a more restrictive definition; they are given by a codimension one foliation $\widetilde{\omega}_\varepsilon$ 
in an infinitesimal neighborhood of $\PP^n$, such that its restriction to the central fiber gives the original form $\omega$ up to scalar multiplication.
As we will see later in Section \ref{unf-def}, if $\omega$ is an integrable global section of $\Omega^1_{\PP^n}(e)$, then first order deformations and unfoldings can be parametrized by the $\CC$-vector spaces
\begin{equation*}
\left.
\begin{aligned}
 D(\omega) &= \left\{\eta\in H^0\left(\Omega^1_{\PP^n}(e)\right):\ \omega\wedge d\eta + d\omega\wedge \eta = 0\right\}\big/\CC.\omega \\
U(\omega) &= \left\{(h,\eta)\in H^0\left(\left(\OO_{\PP^n}\times \Omega^1_{\PP^n}\right)(e)\right):\ hd\omega  = \omega\wedge (\eta - dh) \right\}\big/\CC.(0,\omega)
\end{aligned}
\right.
\end{equation*}
respectively.
As expected, the vector space $D(\omega)$ can be identified with the Zariski tangent space at $\omega$, see \cite{fji} or \cite{hart}.  Denoting by $K(\omega)$ the integrating factors of $\omega$ as in Definition \ref{int-fact}, both types of perturbations can be related {\it via} the exact sequence
%\cite[Section 2.1, p.~691]{fji}
\begin{equation*}
\xymatrix@R=0pt{
0\ar[r] & K(\omega) \ar[r]^-{} & U(\omega) \ar[r]^{} & D(\omega).
}
\end{equation*}

C. Camacho and A. Lins-Neto, in  \cite{caln}, define the following notion of regularity of an integrable, homogeneous, differential 1-form and prove an associated stability result. 
By looking at $\omega$ as a homogeneous affine form in $\CC^{n+1}$, $\omega$ is said to be regular if for every $a<e$ the graded complex of homogeneous elements
\begin{equation}\label{reg-intro}
\xymatrix@R=0pt{
T_{\CC^{n+1}}(a-e) \ar[r]^-{} & \Omega^1_{\CC^{n+1}}(a) \ar[r]^-{} & \Omega^2_{\CC^{n+1}}(a+e)\\
X \ar@{|->}[r] & L_X(\omega) & &\\
& \save[] *{\eta\ }="s" \restore & \save[] +<40pt,0pt> *{\ \omega\vartriangle\eta:=\omega\wedge d\eta + d\omega\wedge \eta}="t" \restore & \ar@{|->}"s";"t"
}
\end{equation}
has trivial homology in degree 1, where $L_X(\omega)$ is the Lie derivative of $\omega$ with respect to the vector field $X$ and in parenthesis we indicate the homogeneous component of the given degree.

\

As a first step towards classifying the space of foliations, we studied the function $a\longmapsto \varphi_{\omega}(a):=dim_\CC\left(Ker(\omega\vartriangle - )(a)\right)$, for every $a\in\NN$. Even if the application $\eta\mapsto\omega\vartriangle \eta$ is a differential operator, we prove in Theorem \ref{triangle-dw} that the values of $\varphi_{\omega}$ grow polynomially. Then, there is a Hilbert polynomial associated to $\varphi_{\omega}$  exposing discrete invariants. The study of the behavior of $\varphi_{\omega}$ and the information attached to its Hilbert polynomial took us to a deeper knowledge of first order unfoldings and deformations that we present in this paper.

\[
\mbox{\LARGE $\star$}
\]

Rational and logarithmic foliations define irreducible components of the space of codimension one foliations, as it is shown by X. G\'omez Mont and A. Lins-Neto in \cite{gmln} and later by F. Cukierman, J. V. Pereira and I. Vainsencher in \cite{fji} for rational foliations and O. Calvo Andrade in \cite{omegar} for logarithmic foliations. Rational and logarithmic foliations in $\PP^n$ can be given, respectively, by differential forms of the type
\[
 \omega_\RR = r F dG - s G dF\qquad\text{and}\qquad\omega_{\LL} = \left(\prod_{i=1}^s f_i \right)\sum_{i=1}^s\lambda_i \frac{df_i}{f_i},
\]
for polynomials $F$ and $G$ of degree $r$ and $s$, and polynomials $f_i$ of degree $d_i$ and scalars $\lambda_i$ such that $\sum_{i=1}^s \lambda_i d_i = 0$ for $s>2$. 

\

In \cite{fji} the authors proves the infinitesimal stability of a generic rational foliation $\omega_\RR$, showing that $D(\omega_\RR)$ is generated by perturbations of the parameters $F$ and $G$. On the other side, using \cite{omegar} one can only show that the perturbations of the parameters $f_i$ and $\{\lambda_i\}$ generate all the deformations of the space of logarithmic foliations as a set, {\it i.e.}, disregarding the sheaf structure along with any algebraic multiplicity that may arise from the equation $\omega\wedge d\omega=0$. This way, we get a partial description of $D(\omega_\LL)$ in terms of the parameters defining $\omega_\LL$, as we will see in Corollary \ref{tangent-logarithmic}.

\

Let us denote by $(\CC^{n+1},0)$ the infinitesimal analytic neighborhood of the origin in $\CC^{n+1}$. Since $\omega$ defines a global foliation in $\PP^n$, it is natural to look at the germ of analytic foliation induced by $\omega$ in $(\CC^{n+1},0)$. The space of analytic germs $U_{hol}(\omega)$ of first order unfoldings of $\omega$, has an analogous definition to the projective one.
As we will see in Section \ref{hol-aff}, one can associate to $U_{hol}(\omega)$ an ideal $I_{hol}(\omega)$ of the space of holomorphic function germs, who allows us to study $U_{hol}(\omega)$ with a nicer structure. We refer to \cite{suwa-review} for a complete exposition on this.

\

In \cite{suwa-meromorphic} and \cite{suwa-multiform}, T. Suwa was able to find generators of $I_{hol}(\omega)$ for generic foliations of rational and logarithmic type, showing that
\begin{equation}\label{suwa-intro}
 I_{hol}(\omega_\RR) = (F,G)\qquad \text{and}\qquad I_{hol}(\omega_\LL) = (F_1,\ldots,F_s)
\end{equation}
where $F_i = \prod_{j\neq i} f_j$.
The lack of this ideal structure in the global projective case, makes difficult to recover the information made available by eq. \ref{suwa-intro} to study first order unfoldings of rational and logarithmic foliations.

To bypass this situation, we propose a new definition of a graded module $\UU(\omega)$ over the ring $S$ of homogeneous coordinates in $\PP^n$, that we call the module of graded projective unfoldings associated to $\omega$. If $R$ denotes the radial vector field and $\Omega^1_S$ the module of differential 1-forms of $S$ over $\CC$,  we define $\UU(\omega)$ as
\[
\UU(\omega) = \left\{(h,\eta)\in S\times \Omega^1_{S}: \ L_R(h)\ d\omega = L_R(\omega)\wedge(\eta - dh) \right\}\big/ S.(0,\omega),
\]
and refer to Definition \ref{action} for the details on the module structure of $\UU(\omega)$. With this module, we can emulate the situation in the infinitesimal analytic case and define, by projecting on the first coordinate, a graded ideal $I(\omega)$ on $S$. By doing so, we can rapidly translate the results of T. Suwa shown in eq. \ref{suwa-intro}, to global foliations in $\PP^n$.

\

As a first application of $\UU(\omega)$, we were able to classify which first order deformations come from first order unfoldings. See Theorem \ref{teo-1} and Theorem \ref{teo-2} for a complete statement of the following results.

\teointro\label{teo1-intro} Let $\omega_\RR$ be a generic rational foliation in $\PP^n$. Then, the following sequence is exact
\begin{equation*}
\xymatrix{
0\ar[r] & K(\omega_\RR) \ar[r]^-{} & U(\omega_\RR) \ar[r]^{} & D(\omega_\RR) \ar[r] & 0.
}
\end{equation*}

\teointro Let $\omega_\LL$ be a generic logarithmic foliation in $\PP^n$ and write as $D(\omega,\overline{f})$ the subspace of $D(\omega_\LL)$ consisting of the perturbations of the parameters $\{f_i\}$. Then, $\pi_2$ is not en epimorphism and the following sequence is exact
\[
\xymatrix@R=0pt{
0\ar[r] & K(\omega_\LL) \ar[r]^-{} & U(\omega_\LL) \ar[r]^{\pi_2} & D(\omega_\LL) \ar[r]^-{} & D(\omega_\LL)/D(\omega_\LL,\overline{f}) \ar[r]^{} & 0.
}
\]

\

It is well known that the singular set of an integrable form $\omega\in H^0(\Omega^1_{\PP^n}(e))$, {\it i.e.} the space $Sing(\omega)$ where $\omega$ vanishes, has always a   codimension two component where $d\omega$ is generically not null, named the Kupka component. The characterization of other components of $Sing(\omega)$ is unknown.

The singular set of foliations of logarithmic type, under some genericity conditions, was studied by F. Cukierman, M. Soares and I. Vainsencher in \cite{fmi} . The authors show that the singular set decomposes as the disjoint union $Sing(\omega)=Z\cup Q$, where $Z$ is the Kupka set of codimension 2 and $Q$ is a finite set of $N_{\omega_\LL}$ points counted with multiplicity. They also give a closed formula to compute $N_{\omega_\LL}$; such computation is done by using a particular formula to obtain the Segre class of the singular scheme of a normal crossing divisor.

\

In Section \ref{Isolated points of the singular set} we give another application of the graded module $\UU(\omega)$ by counting, in a rather simple way, the isolated points of the singular set of logarithmic foliations. See Theorem \ref{teo-3} for a complete statement of:

\teointro\label{teo3-intro} Let $\omega_\LL$ be a generic logarithmic foliation in $\PP^n$ and, as above, decompose $Sing(\omega_\LL) =  Z\cup Q$. Let $\overline{\UU}(\omega)$ denote the classes of isomorphisms of graded projective unfoldings of $\omega$, and $P_{\overline{\UU}(\omega_\LL)}$ its Hilbert polynomial. Then $P_{\overline{\UU}(\omega_\LL)}$ is constant and
\[
 P_{\overline{\UU}(\omega_\LL)}\equiv N_{\omega_\LL} .
\]

\

As always, if $\omega$ defines a codimension 1 foliation in $\PP^n$, in Section \ref{unf-def} we associate to $\omega$ a graded $S$-linear complex $L^\bullet(\omega)$ defined as
\begin{equation*}
\xymatrix@C=30pt{
L^\bullet(\omega): & T_S \ar[r]^-{d\omega\wedge} & \Omega^1_{S} \ar[r]^-{d\omega\wedge} & \Omega^3_{S} \ar[r]^-{d\omega\wedge} & \ldots &
}
\end{equation*}
The module $\UU(\omega)$ and the complex $L^\bullet(\omega)$ are closely related by the following result. See Theorem \ref{phi} for a complete statement:
\teointro Let $\omega$ be a codimension one foliation in $\PP^n$. If we denote by $\mathcal{Z}^1(L^\bullet(\omega))$ and $H^1(L^\bullet(\omega))$ the cycles and homology of $L^\bullet(\omega)$ in degree 1 respectively, then
\begin{equation}\label{teo4-intro} 
 \mathcal{Z}^1(L^\bullet(\omega))/S.\omega \simeq \UU(\omega)\qquad\text{and}\qquad H^1(L^\bullet(\omega))\simeq\overline{\UU}(\omega).
\end{equation}

\

In Section \ref{regularity} we review the notion of regularity of C. Camacho and A. Lins-Neto, defined {\it via} differential operators in terms of the linear complex $L^\bullet(\omega)$.  As an application of $L^\bullet(\omega)$ and the isomorphisms in eq. \ref{teo4-intro} above, we can completely reformulate regularity of rational and logarithmic foliations in terms of unfoldings in the following way. See Theorem \ref{teo-4} for a complete statement:
\teointro Let $\omega$ be a generic rational or logarithmic foliation in $\PP^n$. Then $\omega$ is regular if and only if every first order unfolding of $\omega$ is trivial up to isomorphism.

\

In Section \ref{Isolated points of the singular set} we make some explicit computations with the help of a computer, to show the behavior of $\overline{\UU}(\omega)$ on low degrees. We finally add an appendix where we give a simplified proof of Theorem \ref{teo1-intro}.

\subsection{Acknowledgements}

\hspace{9pt}
I would like to give special thanks to my advisor Fernando Cukierman for the very useful discussions, ideas and suggestions.
I also want to thank Federico Quallbrunn, Cesar Massri, Manuel Dubinsky, Matias del Hoyo, Alicia Dickenstein and Tatsuo Suwa.

The content  of this work is part of the author's doctoral thesis at {\it Universidad de Buenos Aires} under the advice of Fernando Cukierman.

% \tableofcontents
% \phantomsection

\section{Preliminaries}
\label{preliminaries}

\hspace{9pt}Along this section we present the basic definitions and results that we are going to use on the rest of the work. 

\subsection{Codimension 1 algebraic foliations}

\hspace{9pt}In the definition of foliations that we propose, we stress the relative nature of differential forms. The reason for doing so, is for being able to distinguish first order deformations and unfoldings just as foliations over different base spaces.

\

Let us consider two algebraic varieties $T$ and $B$, such that $T$ is of finite type over $B$. We will write $\Omega^1_{T|B}(\LL)$ for the twisted sheaf of differentials 1-forms of $T$ over $B$, for some invertible sheaf $\LL$ on $T$.

\deff We will say that a generically rank 1 subsheaf $\FF=(\omega)$ of $\Omega^1_{T|B}(\LL)$ is an {\it algebraic foliation of codimension 1} on $T$ over $B$, if $\omega$ is a non zero global section of $H^0(\Omega^1_{T|B}(\LL))$ generating $\FF$, such that verifies the Frobenius integrability condition $\omega\wedge d\omega = 0$. We will write $\FF^1(T|B)(\LL)$ for the space of this foliations.

\

In the case where $B=Spec(\CC)$, we will just write $\Omega^1_T(\LL)$ and $\FF^1(T)(\LL)$.
Given a foliation $\FF=(\omega)$, any multiple of $\omega$ by a global section of $\OO_T^*$ defines the same foliation. We then have
\[
\FF^1(T|B)(\LL) = \{\omega\in H^0(\Omega^1_{T|B}(\LL))\big/H^0(\OO_T^*):\ \omega\wedge d\omega = 0 \}.
\]

We are primarily interested in the case where $T=\PP^n$, in this case a foliation $\FF=(\omega)$ is given by a subsheaf of $\Omega^1_{\PP^n}(e):=\Omega^1_{\PP^n}\otimes\OO_{\PP^n}(e)$, for some $e\geq 2$. For such a foliation, the {\it degree} is defined to be the number of common  tangencies with a generic line in $\PP^n$ which is equal to $e-2$. 

Fixing homogeneous coordinates $x_0,\ldots,x_n$, let us denote $S=\CC[x_0,\ldots,x_n]$ the {\it ring of homogeneous coordinates} of $\PP^n$ and $\Omega^r_S$ the module of {\it r-differential forms} of $S$ over $\CC$. The global section $\omega$ can be written as
\[
 \begin{aligned}
  \omega=\sum_{i=0}^n A_i dx_i\in\Omega^1_S
 \end{aligned}
\]
where the $A_i$'s are homogeneous polynomials of degree $e-1$, that verify the integrability condition $\omega\wedge d\omega=0$ and the property of descent to projective space. This condition can be stated as the vanishing of the contraction of $\omega$ with the radial field $R=\sum_{i=0}^nx_i\frac{\partial}{\partial x_i}$.

As we are going to fix one generator for each foliation we might refer simply as $\omega$ to the foliation $\FF=(\omega)$ and note the space of codimension 1 foliations of degree $e-2$ in $\PP^n$ as $\FF^1(\PP^n)(e)$.

\

The quotient $H^0\left(\Omega^1_{\PP^n}(e)\right)\big/H^0\left(\OO^*_{\PP^n}\right)$ identifies with $\PP^N$, for a suitable $N$, by looking at the scalar coefficients  of a differential 1-form $\omega$. Then, the equation $\omega\wedge d\omega = 0$ defines an homogeneous ideal $F$ in such coefficients. The algebraic variety structure of $\FF^1(\PP^n)(e)$ is then given by $\text{Proj}(S/F)\subset \PP^N$.

\

 The {\it Koszul complex} associated to $\omega\in\Omega^1_S$, $K^\bullet(\omega)$, can be defined as
\begin{equation*}
\xymatrix@C=30pt{
K^\bullet(\omega): & S \ar[r]^-{\omega\wedge} & \Omega^1_{S} \ar[r]^-{\omega\wedge} & \Omega^2_{S} \ar[r]^-{\omega\wedge} & \ldots &
}
\end{equation*}
We can use $K^\bullet(\omega)$ to compute the codimension of the singular set of $\omega$, by the well known result:
\teo For $\omega\in \Omega^1_S$ the following are equivalent:
\begin{enumerate}
\item $codim(Sing(\omega))\geq k$
\item $H^l(K^\bullet(\omega))= 0$ for all $l<k$
\end{enumerate}
\begin{proof}
 See \cite[Appendix, p.~172]{m-f1} or \cite[Appendix, p.~87]{m-f2} for two proofs with different level of generalities in the local holomorphic setting and \cite[Theorem 17.4, p.~424]{eisenbud} for a purely algebraic proof of our statement. %and Proposition 18.4, p.~450]{eisenbud} 
\end{proof}

\rem If $\omega\in\FF^1(\PP^n)(e)$ we always have $H^2(K^\bullet(\omega))\neq 0$. This can be seen by looking at the class of $d\omega$ in $H^2(K^\bullet(\omega))$; the integrability condition on $\omega$ makes $d\omega$ a 2-cycle and, by a matter of degree in the homogeneous polynomial coefficients, it can not be border. 

The homology in degree 1 is trivial only in the case where $\omega$ is {\it irreducible}, {\it i.e.}, if it is not of the form $f.\omega'$, for some not invertible function $f$ and a 1-form $\omega'$.

\subsection{Unfoldings and deformations}
\label{unf-def}

\hspace{9pt}Let us write $\CC[\varepsilon]$ for the ring of {\it dual numbers} $\CC[t]/(t^2)$ and $D=Spec(\CC[\varepsilon])$ for the {\it infinitesimal neighborhood of order one}. Let us consider the morphism $i:\PP^n\to\PP^n\times D$, defined by the inclusion to the closed point of $D$.

\deff A {\it first order deformation} of a foliation $\FF=(\omega)\in\FF^1(\PP^n)(e)$, is given by a foliation $\FF_\varepsilon = (\omega_\varepsilon)\in\FF^1(\PP^n\times D|D)(\OO_{\PP^n\times D}(e))$ such that $i^*\FF_\varepsilon\simeq\FF$. We can synthesize this situation with the commutative diagram
\begin{equation*}
\label{pullback1-def}
\xymatrix{
\save[] +<0pt,15pt> *{\FF\ }="t" \restore\PP^n \ar[r]^-{i} \ar[d]_{} & \PP^n\times D\ar[d]^{\pi} \save[] +<0pt,15pt>  *{\ \FF_\varepsilon}="s" \restore \ar@{|-->}_-{i^*} "s";"t"\\ 
Spec(\CC) \ar[r]_{} & D
}
\end{equation*}
where $\pi$ is the projection to $D$. We will say deformation or first order deformation indistinctly.

\

The condition $i^*\FF_\varepsilon\simeq \FF$ allows us to choose $\omega_\varepsilon = \omega +  \varepsilon \eta$ where $\eta\in H^0(\Omega^1_{\PP^n}(e))$. The integrability condition applied to $\omega_\varepsilon$, gives the formula
\[
\omega\vartriangle\eta = 0
\]
recalling that we write $\omega\vartriangle\eta$ for $\omega\wedge d\eta + d\omega\wedge \eta$, as in eq. \ref{reg-intro}.

\

\deff A {\it first order unfolding} of a foliation $\FF=(\omega)\in\FF^1(\PP^n)(e)$ is given by a foliation $\widetilde{\FF}_\varepsilon = (\widetilde{\omega}_\varepsilon)\in\FF^1(\PP^n\times D)(\OO_{\PP^n\times D}(e))$ such that $i^*\widetilde{\FF}_\varepsilon\simeq\FF$. We can synthesize this situation with the commutative diagram
\begin{equation*}
\xymatrix{
\save[] +<0pt,15pt> *{\FF\ }="t" \restore\PP^n \ar[r]^-{i} \ar[d]_{} & \PP^n\times D\ar[d]^{} \save[] +<0pt,15pt>  *{\ \widetilde{\FF}_\varepsilon}="s" \restore \ar@{|-->}_-{i^*} "s";"t"\\  
Spec(\CC) \ar[r]_{} & Spec(\CC)
}
\end{equation*}
We will say unfolding or first order unfolding indistinctly.

\

The condition $i^*\FF_\varepsilon\simeq \FF$ allows us to choose $\widetilde{\omega}_\varepsilon = \omega +  \varepsilon \eta + h d\varepsilon$ where $\eta\in H^0(\Omega^1_{\PP^n}(e))$ and $h\in H^0(\OO_{\PP^n}(e))$. The integrability condition applied to $\widetilde{\omega}_\varepsilon$ can be computed as
\[
\begin{aligned}
\widetilde{\omega}_{\varepsilon}\wedge d\widetilde{\omega}_{\varepsilon} &= (\omega + \varepsilon\eta + hd\varepsilon)\wedge d(\omega +  \varepsilon\eta+ hd\varepsilon) = \\
&= \omega\wedge d\omega + \varepsilon (\omega\vartriangle\eta)+ (\ hd\omega - \omega\wedge(\eta-dh)\ )\wedge  d\varepsilon = 0 .
\end{aligned}
\]

We then have
\begin{equation*}
\widetilde{\omega}_{\varepsilon}\wedge d\widetilde{\omega}_{\varepsilon}=0 \iff 
\left\{
\begin{aligned}
&\omega\vartriangle \eta = 0\\
&hd\omega  = \omega\wedge (\eta - dh) .
\end{aligned}
\right.
\end{equation*}

\prop \label{equiv-equ}Following the situation above, if $\ hd\omega  = \omega\wedge (\eta - dh)$ then $\omega\vartriangle \eta = 0$.
\begin{proof}
If we apply the exterior differential to $hd\omega  = \omega\wedge (\eta - dh)$ we get $2dh\wedge d\omega = -\omega\wedge d\eta + d\omega\wedge\eta$.
Instead, if we multiply $hd\omega  = \omega\wedge (\eta - dh)$ by $\eta-dh$ we get $dh\wedge d\omega = d\omega\wedge\eta$.
Putting together both formulas we find our result.
\end{proof}

In both cases, perturbing in the direction of $\omega$, {\it i.e.} taking $\eta=\omega$, defines the trivial deformation or unfolding.

\deff Let $\FF=(\omega)\in\FF^1(\PP^n)(e)$. We define the $\CC$-vector spaces {\it parameterizing deformations and unfoldings} as
\[
\begin{aligned}
D(\omega) &= \left\{\eta\in H^0\left(\Omega^1_{\PP^n}(e)\right):\ \omega\vartriangle\eta = 0\right\}\big/\CC.\omega\\
U(\omega) &= \left\{(h,\eta)\in H^0\left(\left(\OO_{\PP^n}\times \Omega^1_{\PP^n}\right)(e)\right):\ hd\omega  = \omega\wedge (\eta - dh)\right\}\big/\CC.(0,\omega).
\end{aligned}
\]
We will call $D(\omega)$ and $U(\omega)$ deformations and unfoldings respectively, if no confusion can arise.

\deff Two deformations (unfoldings) $\FF_\varepsilon$ and $\FF_\varepsilon'$ ($\widetilde{\FF}_\varepsilon$ and $\widetilde{\FF}_\varepsilon'$) of a given foliation $\FF\in\FF^1(\PP^n)(e)$ are said to be {\it isomorphic}, if there is an isomorphism $\phi:\PP^n\times D\to \PP^n\times D$ such that
\begin{equation}\label{iso-def-un}
i^*\phi = Id_{\PP^n} \hspace{.5cm} \qquad \text{and}\qquad \FF_\varepsilon\simeq \phi^*\FF_\varepsilon' \hspace{.5cm}\left(\widetilde{\FF}_\varepsilon\simeq \phi^*\widetilde{\FF}_\varepsilon'\right).
\end{equation}

The structure sheaf $\OO_{\PP^n\times D}$ is isomorphic to two copies of $\OO_{\PP^n}$, identifying $1$ and $\varepsilon$ with the canonical vectors. This way, we can decompose an  isomorphism $\phi$ of deformations or unfoldings, as $\phi = \phi_1 + \varepsilon\phi_2$, for $\phi_1,\phi_2\in PGL(n,\CC)$.

The first condition in eq. \ref{iso-def-un} allows us to write $\phi = Id_{\PP^n} + \varepsilon \phi_2$.
Let us now consider the vector field $X=\sum_{i,j} \left(\phi_2\right)_{ij} x_i\frac{\partial }{\partial x_j}$ induced by $\phi_2$, and write $i_X$ for the contraction with $X$. If  $\omega_\varepsilon = \omega + \varepsilon\eta$ and $\widetilde{\omega}_\varepsilon = \omega + \varepsilon\eta + hd\varepsilon$ defines a deformation and an unfolding of $\omega\in\FF^1(\PP^n)(e)$ respectively, by straight forward computation, we get the formulas
\begin{enumerate}
\item $\phi^*\omega_\varepsilon = \omega + \varepsilon (L_X(\omega) + \eta)$
\item $\phi^*\widetilde{\omega}_\varepsilon = \omega + \varepsilon (L_X(\omega) + \eta) + (i_X\omega + h)d\varepsilon$.
\end{enumerate}

\

Denote by $T_{\PP^n}$ the {\it tangent sheaf} in $\PP^n$. 

\deff The spaces of {\it deformations and unfoldings  modulo isomorphism} of $\FF=(\omega)\in\FF^1(\PP^n)(e)$  are the quotients $\overline{D}(\omega):= D(\omega)/C_D(\omega)$ and $\overline{U}(\omega):= U(\omega)/C_U(\omega)$
where
\[
\begin{aligned}\label{clas-def}
C_D(\omega) &= \left\{ L_X(\omega): \ X\in H^0\left(T_{\PP^n}(0)\right)\right\} \\
C_U(\omega) &= \left\{ (i_X\left(\omega),L_X(\omega)\right): \ X\in H^0\left(T_{\PP^n}(0)\right)\right\}.
\end{aligned}
\]

\

\deff\label{int-fact} The space $K(\omega)$ of {\it integrating factors} of $\FF=(\omega)\in\FF^1(\PP^n)(e)$ is given by 
\[
 K(\omega) = \left\{ F\in H^0\left(\OO_{\PP^n}(e)\right) :\  Fd\omega = \omega\wedge \left(-dF\right)\right\} .
\]

\

By Definition \ref{int-fact} and Proposition \ref{equiv-equ}, we  immediately see the exactness of the following sequence, that relates integrating factors, unfoldings and deformations of $\omega$:
\begin{equation}\label{suc_0}
\xymatrix@R=0pt{
0\ar[r] & K(\omega) \ar[r]^-{i_1} & U(\omega) \ar[r]^{\pi_2} & D(\omega) \\
& h \ar@{|->}[r] & (h,0) &\\
& & (h,\eta) \ar@{|->}[r] &\eta
}
\end{equation}

\subsection{Local setting and Cartan's Magic Formula}
\label{hol-aff}

\hspace{9pt}Given a foliation defined by $\omega\in\FF^1(\PP^n)(e)$, we can look at the local foliation induced by $\omega$ by restricting to an open set of $\PP^n$, or by pullbacking $\omega$ to the affine cone in $\CC^{n+1}$. Adopting the former procedure, we keep the homogeneity of $\omega$ and we are able to grade the spaces of unfoldings and deformations in the local algebraic setting in $\CC^{n+1}$, or in the holomorphic infinitesimal setting in $(\CC^{n+1},0)$. We show this in Section \ref{local-setting} below and fix some notation. Then, in Section \ref{cartan}, we recall Cartan's Magic Formula. With this formula, we can decompose affine differential forms as a closed form plus a form which descends to projective space. This last decomposition, which is elementary, is crucial to linearize the notion of unfolding and to connect unfoldings with the notion of regularity, as we will do in Section \ref{sec-3-2} and Section \ref{regularity}, respectively.

\subsubsection{Local setting}\label{local-setting}
\hspace{9pt}We will denote with a subscript $hol$ the analogous definitions  with the preceding section of {\it deformations}, {\it unfoldings} and {\it isomorphism classes of unfoldings}, relative to the space of {\it germs of holomorphic foliations} in $(\CC^{n+1},0)$.

For example, for a germ of holomorphic foliation $\upsilon$ in $(\CC^{n+1},0)$, we have
\[
\begin{aligned}
U_{hol}(\upsilon) &= \left\{(h,\eta)\in \OO_{(\CC^{n+1},0)}\times\Omega^1_{(\CC^{n+1},0)}: hd\upsilon  = \upsilon\wedge (\eta - dh)\right\}\big/\OO_{(\CC^{n+1},0)}.(0,\upsilon)
\end{aligned}
\]
where $\OO_{(\CC^{n+1},0)}$ and $\Omega^1_{(\CC^{n+1},0)}$ are, respectively, the {\it germs of holomorphic functions} and differential {\it 1-forms} in $(\CC^{n+1},0)$.
By projecting the first coordinate of $U_{hol}(\upsilon)$, we get the ideal $I_{hol}(\upsilon)\subset\OO_{(\CC^{n+1},0)}$
\begin{equation}\label{I-J-hol}
\begin{aligned}
 I_{hol}(\upsilon) &= \left\{ h\in \OO_{(\CC^{n+1},0)}: \ hd\upsilon = \upsilon\wedge\widetilde{\eta}\text{ for some }\widetilde{\eta}\in\Omega^1_{(\CC^{n+1},0)}\ \right\}.
%  J_{hol}(\upsilon) &= \left\{ i_X\widetilde{\eta}\in \OO_{(\CC^{n+1},0)}:\ \widetilde{\eta}\in\Omega^1_{(\CC^{n+1},0)}, \  X\in T_{(\CC^{n+1},0)}\ \right\}.
\end{aligned}
 \end{equation}

As we will prove later in Proposition \ref{I/J}, in the case where $\upsilon$ is an irreducible foliation, there is an isomorphism $U_{hol}(\upsilon)\simeq I_{hol}(\upsilon)$. 

\

Let us denote $T_S$ to the module of {\it vector fields}, which we define as $T_S=\Omega^{-1}_S$ the dual of the module of differential 1-forms. We assign to the elements $dx_i$ and $\frac{\partial}{\partial x_i}$ degrees $+1$ and $-1$ respectively.
As in the introduction, if $M$ is a graded $S$-module or a $\CC$-vector space, we will write $M(k)$ for its {\it homogeneous component} of degree $k$. 

By changing in the preceding section $\PP^n$ by $\CC^{n+1}$, we end up with the definitions of {\it algebraic deformations}, {\it unfoldings} and {\it isomorphism classes of unfoldings}. We will note them with a subscript $alg$. 

\

By pullbacking with the application $\pi:\CC^{n+1}\backslash\{0\}\rightarrow\PP^n$, a differential form $\omega\in\FF^1(\PP^n)(e)$ defines a foliation in the affine space $\CC^{n+1}$, or in $(\CC^{n+1},0)$. We will commit a small abuse of notation and keep writing as $\omega$ the pullbacked differential form, seen in $\CC^{n+1}$ or in $(\CC^{n+1},0)$ as well. Being $\omega$ homogeneous, the spaces of holomorphic and algebraic unfoldings are enriched with a natural graded structure. Then, we can decompose
\begin{equation}\label{dec-1}
U_{hol}(\omega) = \prod_{a\in\NN}U_{hol}(\omega)(a)\qquad \text{and}\qquad U_{alg}(\omega) = \bigoplus_{a\in\NN}U_{alg}(\omega)(a),
\end{equation}
where $U_{hol}(\omega)(a)$ and $U_{alg}(\omega)(a)$ can be readily identified. 

\remo Taking $a=e$, the degree of $\omega$, we have isomorphisms
\begin{equation}\label{dec-2}
U_{hol}(\omega)(e)\simeq U_{alg}(\omega)(e) \simeq U(\omega)
\end{equation}
showing that any unfolding of a projective foliation in $\PP^n$, can be obtained as the homogeneous component of some holomorphic germ of unfolding in $(\CC^{n+1},0)$, or, also, as some algebraic unfolding in $\CC^{n+1}$. 
Proceeding in an analogous way to equations (\ref{dec-1}) and (\ref{dec-2}), we can conclude that the same statement holds for algebraic deformations in $\PP^n$.

\subsubsection{Cartan's Magic Formula}\label{cartan}

\hspace{9pt}Following {\it Cartan's  Magic Formula}, we can compute compute the Lie derivative of a differential form $\tau$, with respect to a vector field $X$, as
\begin{equation*}\label{cartan-formula0}
 L_X(\tau) = i_Xd\tau + di_X\tau
\end{equation*}
see {\it e.g}. \cite{warner}.

Let us take $\tau\in\Omega^r_S(p)$ and $R\in T_S(0)$, the radial vector field. In this case, the formula above gives the equality 
\begin{equation}\label{lie-desc}
L_R(\tau) = di_R\tau+i_R d\tau = p\tau
\end{equation}
which allows us to decompose $\tau=\tau_d + \tau_r$, where $\tau_d$ and $\tau_r$ are the exact and radial terms respectively, see {\it e.g.} \cite{jou}.

\

Let us recall that if $\FF$ is an $\OO_{\PP^n}$-module, then the functor $\Gamma_*$ defines a graded $S$-module as $\Gamma_*(\FF)=\bigoplus_{a\in\ZZ}H^0(\FF(a))$. By looking at eq. \ref{lie-desc}, we can define the graded morphism $\CCal:\Omega^1_S\to\Gamma_*(\OO_{\PP^n}\times\Omega^1_{\PP^n})$ as $\CCal=\bigoplus_{a\in\NN}\CCal_a$, where $\CCal_a$ is given by the formula
\begin{equation}\label{cartan-formula}
\xymatrix@R=5pt@C=50pt{
\Omega^1_S(a) \ar[r]^-{\CCal_a} & H^0\left(\left(\OO_{\PP^n}\times\Omega^1_{\PP^n}\right)(a)\right)\\
\save[] *{\eta\ }="s" \restore& \hspace{80pt}\save[]  +<2pt,0pt> *{\ \left(-\frac{1}{a}i_R\eta,\frac{1}{a}i_Rd\eta\right)}="t" \restore \\ \ar@{|->} "s";"t"
}
\end{equation}
The following property is immediate:
\prop The application $\CCal:\Omega^1_S\to\Gamma_*(\OO_{\PP^n}\times\Omega^1_{\PP^n})$ defines an isomorphism in each homogeneous component.

\

We will usually write $\CCal_a(\eta) = (h,\eta_r)$. We adopt the minus sign to be able to write $\eta = \eta_r - dh$ which has a direct relation with the definition of unfolding. 

\subsection{Rational and logarithmic foliations}

\hspace{9pt}Along this section, we present the definitions of rational and logarithmic foliations and we fix some genericity conditions. Then, we recall some important results attached to this type of foliations that we will use later: the characterization of the first order deformations and the characterization of the ideal $I_{hol}$ associated to first order unfoldings.

\deff A {\it rational foliation} of type $(r,s)$ in $\FF^1(\PP^n)(e)$, is defined by an $\omega_\RR\in H^0(\Omega^1_{\PP^n}(e))$ of the form
\begin{equation*}\label{rational}
\omega_{\RR} = rFdG- sGdF,
\end{equation*}
where $F$ and $G$ are homogeneous polynomials of degrees $r$ and $s$ respectively, and $r+s=e$. The Zariski clousure in $\FF^1(\PP^n)(e)$ of foliations of this type defines the {\it set of rational foliations} which will be denoted as $\RR(n,(r,s))$. We define the generic open set $\mathcal{U}_\RR\subset \RR(n,(r,s))$ as 
\begin{equation}\label{gen-rac}
\mathcal{U}_\RR = \{\omega\in \RR(n,(r,s)):\ codim(Sing(d\omega))\geq 3,\ codim(Sing(\omega))\geq 2 \}.
\end{equation}

First order deformations of rational foliations are studied in the works \cite{gmln} and \cite{fji}. The latter, takes into account the scheme structure of codimension one foliations and proves, among other things, that $\FF^1(\PP^n)(e)$ is generycally reduced at a rational foliation. We recall from \cite[Proposition 2.4, p.~693]{fji} the following result.
\teo \label{tangent-rational} Let $\omega_{\RR}\in\mathcal{U}_{\RR}\subset\FF^1(\PP^n)(e)$. Then, the first order deformations of $\omega_{\RR}$ are given by the perturbations of the parameters $F$ and $G$
\begin{displaymath}
D(\omega_\RR) = Span\left(\left\{ \eta\in\RR(n,(r,s)):\eta = rfdG- sGdf \text{ or }\eta = rFdg- sgdF \right\}\right)\big/ \CC.\omega_{\RR}.
\end{displaymath}

\

In the case of germs of holomorphic foliations in $(\CC^{n+1},0)$, let us refer to a rational foliation as generic in an analogous sense to eq. \ref{gen-rac}. We can recall from \cite[Proposition 1.7, p.~102]{suwa-meromorphic} the following result.

\teo \label{I-meromorphic}Let $\upsilon\in\Omega^1_{(\CC^{n+1},0)}$  define a generic rational foliation in $(\CC^{n+1},0)$. If $\upsilon$ is of the form $fdg-gdf$, then $I_{hol}(\upsilon) = (f,g)$. 

\

\deff A {\it logarithmic foliation} of type $(d_1,\ldots,d_s)$  in $\FF^1(\PP^n)(e)$, is defined by an $\omega_\LL\in H^0(\Omega^1_{\PP^n}(e))$ of the form 
\begin{equation}\label{logarithmic}
\omega_{\LL} = \left(\prod_{i=1}^s f_i \right)\sum_{i=1}^s\lambda_i \frac{df_i}{f_i},
\end{equation}
where $s\geq 3$ and 

\begin{enumerate}
\item $(\lambda_1,\ldots,\lambda_s)\in\Lambda(s) :=\{(\lambda_1,\ldots,\lambda_s)\in\CC^s:\  \lambda_1d_1+\ldots+\lambda_s d_s=0\}$ 
\item $f_i$ is homogeneous of degree $d_i$ and $d_1+\ldots+d_s=e$.
\end{enumerate}
The Zariski clousure in $\FF^1(\PP^n)(e)$ of foliations of this type defines the {\it set of logarithmic foliations} which will be denoted as $\LL(n,\overline{d})$.  We define the generic open set $\mathcal{U}_\LL\subset \LL(n,\overline{d})$ as
\begin{equation}\label{gen-log}
\mathcal{U}_\LL = \left\{\omega\in \LL(n,(\overline{d})):\ \omega\text{ verifies a) and b) below }\right\},
\end{equation}
writing $\omega=\left(\prod_{i=1}^s f_i \right)\sum_{i=1}^s\lambda_i \frac{df_i}{f_i}$ we have the conditions:
\begin{enumerate}
\item[a)] $D = \{f_1.\ldots . f_s = 0\}$ is a normal crossing divisor
% \item[b)] the height of the ideal $(f_i,f_j,f_k)$ is 3 for every triplet $i,j,k$
\item[b)] $\lambda_i \neq \lambda_j (\neq 0)$ for every $i\neq j$.
\end{enumerate}

\

We will usually note $\overline{d}$, $\overline{\lambda}$ and $\overline{f}$ the $s$-uples involved in the expression of a logarithmic foliation.
Noting $F_i = \prod_{j\neq i} f_j$, we will frequently write $\omega_{\LL}$ as
\begin{equation*}\label{logarithmic2}
\omega_{\LL} = \sum_{i=1}^s \lambda_i\  F_i \ df_i.
\end{equation*}

Let us fix $\omega_{\LL}$ as in eq. \ref{logarithmic} and define the spaces of perturbation of parameters of $\omega_{\LL}$ as

\[
\begin{aligned}
D(\omega_\LL,\overline{f}) &= Span\left(\{\eta_{g_i}\in\LL(n,\overline{\lambda}):\ \eta_{g_i}\text{ equals }\omega_{\LL}\text{ with }f_i\text{ changed by }g_i\}\right)\big/\CC.\omega_\LL\\
D(\omega_\LL,\overline{\lambda}) &= Span\left(\{\eta_{\overline{\mu}}\in\LL(n,\overline{\mu}):\ \eta_{\overline{\mu}}\text{ equals }\omega_{\LL}\text{ with }\overline{\lambda}\text{ changed by }\overline{\mu}\}\right)\big/\CC.\omega_\LL.
\end{aligned} 
\]

By direct computation, it is straight forward to check that $D(\omega_\LL,\overline{f})$ and $D(\omega_\LL,\overline{\lambda})$ are subspaces of $D(\omega_\LL)$.

\

Logarithmic foliations has been studied in \cite{omegar} where it is shown that they are an irreducible component of the space of codimension one foliations. The analytic (set theoretical) approach of this work, allows us to compute first order deformations of $\omega_\LL$ in the space of foliations with its \emph{reduced scheme structure} $\FF^1_{red}(\PP^n)(e)$. 

Regarding the scheme structure of $\FF^1(\PP^n)(e)$, there is an ongoing work by F. Cukierman \emph{et al.}, see \cite{fjc}, where they show that $\FF^1(\PP^n)(e)$ is generycally reduced at a logarithmic foliation. We will not use this result at all. However, the effects of the reduced structure of $\FF^1(\PP^n)(e)$ can be immediately understood in our statements.

\

To take into account Calvo's results in $\FF^1(\PP^n)(e)$, first we give to the set of codimension one foliations in $\PP^n$ the reduced structure $\FF^1_{red}(\PP^n)(e)$, see \cite[II, \S 3, Theorem 2, p.~88]{mumford}, and consider the natural inclusion
%Through this identification, since $\LL(n,\overline{d})$ consists only of schematically closed points, logarithmic foliations are the same set in both settings.
%\cite[II, Proposition 2.6, p.~78]{hart}
\[
\xymatrix{
\FF^1_{red}(\PP^n)(e) \ar@{^{(}->}[r] & \FF^1(\PP^n)(e).
}
\]
This map induces another inclusion of Zariski tangent spaces, see \cite[III, \S 4, pp.~170-171]{mumford},
\[
\xymatrix{
T_{\omega}\FF^1_{red}(\PP^n)(e) \ar@{^{(}->}[r] & T_\omega\FF^1(\PP^n)(e)
}
\]
%  or \cite[Secion 2.1, p.~691]{fji}
for every closed point $\omega$. Then, by \cite[II, Excercise 2.8, p.~80]{hart}, we can always identify first order deformations with Zariski tangent spaces
\[
D_{red}(\omega) = T_{\omega}\FF^1_{red}(\PP^n)(e)\qquad \text{and} \qquad D(\omega)=T_{\omega}\FF^1(\PP^n)(e),
\]
and decompose $D(\omega)$ as
\begin{equation}\label{desc-t}
D(\omega) = D_{red}(\omega)\bigoplus D_{+}(\omega)
\end{equation}
where $D_{+}(\omega):=D(\omega)\big/D_{red}(\omega)$ can be seen as the first order deformations arising from the non-reduced structure of $\FF^1(\PP^n)(e)$ at $\omega$.

\

From \cite{omegar} we can state the following decomposition of $D(\omega_\LL)$.

\teo\label{tangent-logarithmic} Following the notation above, there exists an open Zariski set $\mathcal{U}\subset\LL(n,\overline{d})\subset\FF^1(\PP^n)(e)$ such that if $\omega_\LL\in\mathcal{U}$ then the first order deformations of $\omega_\LL$ can be decomposed as
\[
D(\omega_\LL) = D(\omega_\LL,\overline{f})\bigoplus D(\omega_\LL,\overline{\lambda})\bigoplus D_{+}(\omega_\LL).
\]
\begin{proof}By eq. \ref{desc-t} we just need to show that
\[
D_{red}(\omega_\LL)= D(\omega_\LL,\overline{f})\bigoplus D(\omega_\LL,\overline{\lambda}). 
\]
For doing this, from \cite[Theorem 4.1, p.~766]{omegar} and \cite[Corollary 4.2, p.~766]{omegar}, we can consider the following parametrization map
\[
\xymatrix@R=8pt{
\PP(\Lambda(s))\times \prod\limits_{i=1}^s \PP\left(H^0\left(\OO_{\PP^n}(d_i)\right)\right) \ar[r]^-{\varphi} & \LL(n,\overline{d})\subset \FF^1_{red}(\PP^n)(e)\\
\save[] +<-15pt,-0pt> *{(\overline{\lambda},\overline{f})\ }="s" \restore & \save[] +<-0pt,-0pt> *{\ \sum_{i=1}^s \lambda_i\  F_i \ df_i}="t" \restore  \ar@{|->}"s";"t"
}
\]

The differential of $\varphi$ at a point $(\overline{\lambda},\overline{f})$ can be computed as
\begin{equation}\label{diff-tang}
d\varphi|_{(\overline{\lambda},\overline{f})}(\overline{\mu},\overline{g}) = \eta_{\overline{\mu}} + \sum_{i=1}^s\eta_{g_i}.
\end{equation}
% for a tangent direction $(\overline{\mu},\overline{g})\in\left(\Lambda(s)/\overline{\lambda}\right)\times \left(\prod\limits_{i=1}^s H^0\left(\OO_{\PP^n}(d_i)\right)/(f_i)\right)$, see \cite[Section 2.1. p.~691]{fji}.
This formula its obtained by looking at the pullback 
\[
\xymatrix{
\Lambda(s)\times \prod\limits_{i=1}^s H^0\left(\OO_{\PP^n}(d_i)\right) \ar[r]^-{\pi^*\varphi} & \LL(n,\overline{d})\subset \Omega^1_{S}\big/\CC
}
\]
which is a multilinear application, see \cite[Chap. VIII, 8.1.4 p.~152]{dieudonne}.

Taking bases of the vector spaces involved, the multilinearity of $\pi^*\varphi$ also allows us to express $\pi^*\varphi(\overline{\lambda},\overline{f})$ in terms of algebraic operations on the coordinates of $(\overline{\lambda},\overline{f})$. The same goes for $d\pi^*\varphi$ and, passing to the quotient, to $d\varphi$.  Thus, the determinant of $d\varphi$  will be an isomorphism in an open Zariski set $\mathcal{U}\subset \LL(n,\overline{d})$.

Finally, taking $(\overline{\lambda},\overline{f})\in \varphi^{-1}(\mathcal{U}) $, eq. \ref{diff-tang} shows that $D(\omega_\LL,\overline{f})$ and $D(\omega_\LL,\overline{\mu})$ are in direct sum and, by a dimensional argument, they span all the tangent space of $\FF^1_{red}(\PP^n)(e)$ at $\omega_\LL$.
\end{proof}

% \rem $\LL(n,\overline{d})$ and its {\it reduced structure} $\LL(n,\overline{d})_{red}$ share the same topological spaces, see \cite[Chap. II, Section 5, Proposition 4, p.~110]{mumford}. On the other side, even if the topologies of $\LL_{hol}(n,\overline{d})$ and $\LL(n,\overline{d})_{red}$ are different, by the proof of \cite[Chap. II, Section 2, Theorem 2, p.~88]{mumford}, we can identify open Zariski sets in both spaces.
% 
% \
% 
% Following the above remark we have:
% 
% \cor\label{tangent-logarithmic-cor} Let $\omega_\LL\in\mathcal{U}\subset \LL(n,\overline{d})$ as in Theorem \ref{tangent-logarithmic} above. Then we have the inclusion
% \[
% \xymatrix{
% D(\omega_\LL,\overline{f})\bigoplus D(\omega_\LL,\overline{\lambda}) \ar@{^(->}[r]^{} & D(\omega_\LL).
% }
% \]

\

In the case of germs of holomorphic foliations in $(\CC^{n+1},0)$, let us refer to a logarithmic foliation as generic in an analogous sense to eq. \ref{gen-log}. We can recall from \cite[Proposition 1.7, p. 102]{suwa-multiform} the following result.

\teo \label{I-multiform} Let $\upsilon\in\Omega^1_{(\CC^{n+1},0)}$  define a generic logarithmic foliation in \linebreak$(\CC^{n+1},0)$. If $\upsilon$ is of the form $\upsilon = \sum_{i=1}^s \lambda_i\  F_i \ df_i$, then  $I_{hol}(\upsilon) = (F_1,\ldots,F_s)$.

\section{Graded projective unfoldings}
\label{graded projective unfoldings}

\hspace{9pt}Along this section, let us fix $\omega\in\FF^1(\PP^n)(e)$ and regard it as an affine form in $\Omega^1_S$.

Here we present our main objects of study, which are the module of graded projective unfoldings $\UU(\omega)$ and the linear complex $L^\bullet(\omega)$. The main idea behind $\UU(\omega)$ is to be able to extend local properties to global ones. The complex $L^\bullet(\omega)$ allows us to understand first order unfoldings in terms of a linear operator and to connect them to the notion of regularity.

\subsection{Graded projective unfoldings}

\deff\label{action} We define the $S$-module of {\it graded projective unfoldings} of $\omega$ as
\[
\UU(\omega) = \left\{(h,\eta)\in S\times \Omega^1_{S}: \ L_R(h)\ d\omega = L_R(\omega)\wedge(\eta - dh) \right\}\big/ S.(0,\omega).
\]
For $a\in\NN$, the homogeneous component of degree $a$ can be written as
\[
\UU(\omega)(a) = \left\{(h,\eta)\in (S\times \Omega^1_{S})(a): \ a\ h\ d\omega = e \ \omega\wedge(\eta - dh) \right\}\big/ S(a-e).(0,\omega). 
\]

For $(h,\eta)\in \UU(\omega)(a)$ and $f\in S(b)$, the graded {\it $S$-module structure} is defined {\it via} the formula
\begin{equation*}
f\cdot (h,\eta) := \left( fh,\  \tfrac{(a+b)}{a}\ f\eta \ +\ \tfrac{1}{a}\left(a\ h\ df- b\ f\ dh\right)\ \right)\in\UU(\omega)(a+b).
\end{equation*}

\

\prop\label{radiality} If $(h,\eta)\in\UU(\omega)(a)$, then $(h,\eta)\in H^0\left(\left(\OO_{\PP^n}\times \Omega^1_{\PP^n}\right)(a)\right)$ .
\begin{proof}
By contracting the equation $ah\ d\omega = e\ \omega\wedge(\eta - dh)$ with the radial field $R$, we can see that $i_R\eta=0$. This shows that the pair $(h,\eta)$ defines a global section of $\left(\OO_{\PP^n}\times\Omega^1_{\PP^n}\right)(a)$ as the name of $\UU(\omega)$ suggests.
\end{proof}

\

\deff We define the {\it isomorphism classes} of graded projective unfoldings, as the quotient $\overline{\UU}(\omega):=\UU(\omega)/C_{\UU}(\omega)$. For $a\in\NN$,  an homogeneous component of degree $a$ of $C_{\UU}(\omega)$, is defined as
\[
C_{\UU}(\omega)(a) = \left\{\left(i_X\omega,\frac{a\ i_Xd\omega + e\  di_X\omega}{e}\right): \ X\in T_S(a-e) \right\}\Big/ S(a-e).(0,\omega).
\]

\

Emulating the situation in $(\CC^{n+1},0)$ of eq. \ref{I-J-hol}, we define:

\deff Let $\pi_1:\UU(\omega)\to S$ be the projection to the first coordinate. We define the graded ideals of $S$ associated to $\omega$ as 
\[
 \begin{aligned}
  I(\omega) &= \pi_1(\UU(\omega)) = \left\{ h\in S:\  hd\omega = \omega\wedge\widetilde{\eta}\text{ for some }\widetilde{\eta}\in\Omega^1_S\right\}\\
  J(\omega) &= \pi_1(C_{\UU}(\omega)) = \left\{ i_X(\omega)\in S:\ X\in T_S \right\}.
 \end{aligned}
\]

\rem\label{obs} From eq. \ref{dec-1}, we can see that $I_{hol}(\omega)$ is a graded ideal, so is generated by polynomials. Then, $I(\omega)\simeq I_{hol}(\omega)$ and from Theorem \ref{I-meromorphic} and Theorem \ref{I-multiform} we get the generators of $I(\omega)$ as well.

\prop\label{I/J} The projection $\pi_1:\UU(\omega)\xrightarrow{} S$ induces the isomorphism
\[
\overline{\UU}(\omega)\simeq I(\omega)/J(\omega).
 \]
And, in the case where $\omega$ is irreducible, we also have $\UU(\omega)\simeq I(\omega)$.
\begin{proof} Let us consider $(h,\eta_1),(h,\eta_2)\in (S\times\Omega^1_S)(a)$ such that 
\begin{align*}
a\ hd\omega &= e\ \omega\wedge(\eta_1-dh)\\
a\ hd\omega &= e\ \omega\wedge(\eta_2-dh).
\end{align*}
Then $\omega\wedge (\eta_1-\eta_2) = 0$. In the case where $\omega$ is irreducible, there must exist $f\in S(a-e)$ such that $\eta_1-\eta_2 = f\omega$. This way the classes of $(h,\eta_1)$ and $(h,\eta_2)$ coincide in $\UU(\omega)$, which shows that $\UU(\omega)\simeq I(\omega)$. By doing the same for elements of the form $\left(i_X\omega,\frac{a\ i_Xd\omega + e\  di_X\omega}{e}\right)$ we can see the isomorphism $C_\UU(\omega)\simeq J(\omega)$.

Regardless the irreducibility of $\omega$, putting together both arguments we have that $\overline{\UU}(\omega)\simeq I(\omega)/J(\omega)$.
\end{proof}

\remo For $(h,\eta)\in (S\times\Omega^1_S)(a)$ the application $(h,\eta) \mapsto \left(h,\frac{a\eta + (e-a)dh}{e}\right)$ gives isomorphisms between $U_{alg}(\omega)(a)$ and $\UU(\omega)(a)$. The twisted $S$-module structure of $\UU(\omega)$ is motivated by the ideal structure of $I(\omega)$ seen through this isomorphism.

\subsection{The complex $L^\bullet(\omega)$}\label{sec-3-2}

\hspace{9pt}The equivalence between the conditions $\omega\wedge d\omega = 0$ and $d\omega\wedge d\omega = 0$, allows us to define the following complex:

\deff\label{complex_dw} We define the graded complex  $L^\bullet(\omega)$ of $S$-modules associated to $\omega$, as
\begin{equation*}
\xymatrix@C=30pt{
L^\bullet(\omega): & T_S \ar[r]^-{d\omega\wedge} & \Omega^1_{S} \ar[r]^-{d\omega\wedge} & \Omega^3_{S} \ar[r]^-{d\omega\wedge} & \ldots &
}
\end{equation*}
where $L^s(\omega)=\Omega^{2s-1}_S$ for $s\geq 0$ and the 0-th differential is defined as $d\omega\wedge X:= i_Xd\omega$. 
The grading of $L^\bullet(\omega)$ is given by the decomposition $L^\bullet(\omega) = \bigoplus_{a\in\NN} L^\bullet(\omega)(a)$, where $L^\bullet(\omega)(a)$ is the complex of finite vector spaces
\begin{equation*}
\xymatrix@C=25pt{
L^{\bullet}(\omega,a): & T_S(a-e) \ar[r]^-{d\omega\wedge} & \Omega^1_{S}(a) \ar[r]^-{d\omega\wedge} & \Omega^3_{S}(a+e) \ar[r]^-{d\omega\wedge} & \ldots .
}
\end{equation*}

As usual, we note $\mathcal{Z}^k(-)$ and $\mathcal{B}^k(-)$ for the cycles and borders of degree $k$, respectively, of the given complex.

\

Recall from eq. \ref{cartan-formula} the definition of the isomorphism $\CCal:\Omega^1_S\to \Gamma_*(\OO_{\PP^n}\times\Omega^1_{\PP^n})$ and let us consider the inclusion $i:\Gamma_*(\OO_{\PP^n}\times \Omega^1_{\PP^n})\to S\times\Omega^1_S$.

\teo\label{phi} The composition $i\circ\CCal:\Omega^1_S\to S\times \Omega^1_S$ induces isomorphisms
\[
\mathcal{Z}^1(L^\bullet(\omega))\big/ S.\omega \simeq \UU(\omega)\qquad\text{and}\qquad H^1(L^\bullet(\omega)) \simeq \overline{\UU}(\omega).
\]

\begin{proof}
Let us consider $\eta\in \Omega^1_{S}(a)$ such that $d\omega\wedge\eta = 0$. Applying $i\circ\CCal_a(\eta) = (h,\eta_r)$, we can decompose it as $\eta = \eta_r-dh$. By contracting with the radial field we have
\[
\begin{aligned}
i_R(d\omega\wedge(\eta_r-dh)) &= 0 \iff \\
a\ h\ d\omega &= e\ \omega\wedge(\eta_r - dh).
\end{aligned}
\]

On the other side, consider  a pair $(h,\eta')\in (S\times\Omega^1_S)(a)$ such that
\begin{equation}\label{ecu10}
a\ h\ d\omega = e\ \omega\wedge(\eta'-dh).
\end{equation}
Following Proposition \ref{radiality}, we now that $(h,\eta')$ is in the image of $i\circ\CCal$ for some $\eta=\eta'-dh$. Multiplying eq. \ref{ecu10} by $\eta' - dh$ we obtain $d\omega\wedge(\eta' - dh)=0$.

Since $i\circ\CCal(\omega)=\omega$, passing to the quotient of $S.\omega$ we finally have the isomorphism
\[
\mathcal{Z}^1(L^\bullet(\omega))\big/ S.\omega\simeq \UU(\omega).
\]

Let us consider now an element $\left(i_X\omega,\frac{a\ i_Xd\omega + e\  di_X\omega}{e}\right)\in C_{\UU}(\omega)(a)$. By the equality
\[
\begin{aligned}
\frac{a\ i_Xd\omega + e\  di_X\omega}{e} \ - \ di_X\omega \ =\  \frac{a}{e}\ i_Xd\omega\ = \ d\omega\wedge\left(\frac{a}{e}X\right)
\end{aligned}
\]
We then have that $i\circ\CCal$ also induces an isomorphism between $ \mathcal{B}^1(L^\bullet(\omega))\big/ S.\omega\simeq C_{\UU}(\omega)$ and the result follows.
\end{proof}

\cor\label{H^1I/J} Following the conditions of Theorem \ref{phi}, we can also write $H^1(L^\bullet(\omega))\simeq I(\omega)/J(\omega)$.

\section{Deformations modulo unfoldings}
\label{deformations modulo unfoldigns}
\hspace{9pt}In this section we apply Suwa's local results on rational and logarithmic foliations in the global projective setting using the graded module $\UU(\omega)$. By doing so, we classify which first order deformations arise from first order unfoldings, see Theorems \ref{teo-1} and \ref{teo-2}, respectively.
\subsection{Rational foliations}
\teo\label{teo-1} Let $\omega_\RR\in\mathcal{U}_\RR\subset \RR(n,(r,s))$ be a generic rational foliation. Then, the following sequence is short exact
\begin{equation*}
\xymatrix{
0\ar[r] & K(\omega_\RR) \ar[r]^-{i_1} & U(\omega_\RR) \ar[r]^{\pi_2} & D(\omega_\RR) \ar[r] & 0
}
\end{equation*}

\begin{proof}
By Remark \ref{obs} and  Theorem \ref{I-meromorphic} we have $I(\omega_\RR)=(F,G)$. By the genericity conditions, $\omega_\RR=rFdG-sGdF$ is irreducible and then $\pi_1:\UU(\omega_\RR)\rightarrow I(\omega_\RR)$ is an isomorphism. It is straight forward to check that $\pi_1$ verifies $\pi_1^{-1}(F)=(F,0)\in\UU(\omega_\RR)(r)$ and $\pi_1^{-1}(G)=(G,0)\in\UU(\omega_\RR)(s)$.
 
 Since $U(\omega_\RR)=\UU(\omega_\RR)(e)$, we just need to find which elements appear in $\UU(\omega_\RR)(e)$ by the action of $S$ defined in Definition \ref{action}, applied to $(F,0)$ and $(G,0)$.
 
 Let us consider $g\in S$ of degree $s=e-r$. Multiplying $g\cdot (F,0)$ we obtain
 \[
{g}\cdot (F,0) = \left( {g}F, \ \tfrac{1}{r} \left(\, r\, F\, d{g} \, - \, s\,  {g} \, dF\,\right)\right)\in\UU(\omega_\RR)(e).
\]
In the same way, taking $f\in S(r)$ and multiplying $f\cdot ({G},0)$, we will have
\[
{f}\cdot ({G},0) = \left({f}{G}, \ \tfrac{1}{s} \left(\,s\,{G}\,d{f} \,-\, r\, {f}\,d{G}\,\right)\right)\in\UU(\omega_\RR)(e).
\]
Looking at the second coordinate of this elements and using the classification of $D(\omega_\RR)$ from Theorem \ref{tangent-rational}, the result follows.
\end{proof}

\rem A simplified proof of the above result can be given by checking the equality of the dimensions of the vector spaces $I(\omega_\RR)(e)$ and $D(\omega_\RR)$. Also, we give another proof following the ideas of \cite{suwa-meromorphic}, see Appendix \ref{proof2}, without using Theorem \ref{tangent-rational}. Anyway, we write our previous demonstration not because of its comparison with these two alternative computations, but because of its natural extension to the case of logarithmic foliations, in which case we do not know any other proof to our result.

\subsection{Logarithmic foliations}

From eq. \ref{desc-t} and Theorem \ref{tangent-logarithmic}, recall the open Zariski set $\mathcal{U}$ and the decomposition $D(\omega_\LL)= D(\omega_\LL,\overline{f})\bigoplus D(\omega_\LL,\overline{\lambda})\bigoplus D_+(\omega_\LL).$
% \[
% D(\omega_\LL) = D(\omega_\LL,\overline{f})\bigoplus D(\omega_\LL,\overline{\lambda})\bigoplus D_+(\omega_\LL).
% \]
% 
% \

\teo\label{teo-2} Let $\omega_\LL = \sum_{i=1}^s \lambda_i\  F_i \ df_i\in\mathcal{U}_\LL\cap\mathcal{U}\subset\LL(n,\overline{d})$ be a generic logarithmic foliation. Then $\pi_2$ is not an epimorphism and its image is $D(\omega_\LL,\overline{f})$, making the following sequence to be exact
\[
\xymatrix@R=0pt{
0\ar[r] & K(\omega_\LL) \ar[r]^-{i_1} & U(\omega_\LL) \ar[r]^{\pi_2} & D(\omega_\LL) \ar[r]^-{\pi} & D(\omega_\LL,\overline{\lambda})\bigoplus D_+(\omega_\LL) \ar[r]^{} & 0\\
}
\]
where the last projection is the natural one.

\begin{proof}
From the decomposition of $D(\omega_\LL)$ fo Theorem \ref{tangent-logarithmic}, we just need to show that the image of $\pi_2$ is exactly $D(\omega_\LL,\overline{f})$ to get our result.

By Remark \ref{obs} and Theorem \ref{I-multiform} we have $I(\omega_\LL) = (F_1,\ldots, F_s)$. For $i=1,\ldots,s$ we want to find $\theta_i$ such that $(F_i,\theta_i)\in\UU(\omega_\LL)(b_i)$, where $b_i=e-d_i$. Once we find this elements, we will be able to get the generators of $\UU(\omega_\LL)(e)=U(\omega_\LL)$ using the action of Definition \ref{action}.

One might think that a perturbation induced by $F_i = \prod_{j\neq i}f_j$ is going to be normal to the direction given by $f_i$, and so, that $i_{\frac{\partial}{\partial f_i}}\theta_i =0$. By the transversality of the $\{f_i\}$, we might deal with them as a system of parameters and compute $i_{\frac{\partial}{\partial f_i}}df_j=0$, for $i\neq j$, and $i_{\frac{\partial}{\partial f_i}}df_j=1$, for $i=j$. With this assumptions, fix $i$ and contract the following equation by the vector field $\frac{\partial}{\partial f_i}$:
\[
 b_i\, F_i \,d\omega_{\LL} = e\ \omega_{\LL}\wedge(\theta_i-dF_i)\,,
\]
then, we can effectively clear $\theta_i$ as
\[
 \theta_i = \frac{b_i}{e\lambda_i}\sum_{\substack{j=1\\j\neq i}}^s(\lambda_j-\lambda_i)\ F_{ji} df_j\ + \ dF_i.
\]
Now, it is immediate to see that $\pi_1^{-1}(F_i)=(F_i,\theta_i)\in\UU(\omega_i)(b_i)$.

Let us take $g\in S(d_i)$ and compute the multiplication
\[
g\cdot(F_i ,\theta_i ) = \left(gF_i ,\ \tfrac{e}{b_i}\,g\theta_i  + \tfrac{1}{b_i} \left(\, b_i  \,F_i\, dg\, -\, d_i\, gdF_i  \,\right)\  \right).
\]

Expanding $\theta_i$ in the second coordinate, we found that
\[
\begin{aligned}
&\frac{1}{\lambda_i}g\sum_{\substack{j=1\\j\neq i}}^s\ (\lambda_j-\lambda_i)\ F_{ji}\ df_j \ +\tfrac{e}{b_i}\ {g}\ dF_i \ + \ F_i\ d{g} \ -\  \tfrac{d_i}{b_i}{g}\ dF_i = \\
&\hspace{1.1cm}= \frac{1}{\lambda_i}\left( \sum_{j\neq i}\ \lambda_j\ {g}F_{ji}\ df_j \ + \ \lambda_i\ F_i\ d{g} \right)\in D(\omega_\LL,\overline{f})
\end{aligned}
\]
Even more, $\pi_2(g\cdot(F_i ,\theta_i ))$ is exactly the perturbation of $\omega_\LL$, given by replacing $f_i$ by $g$. Doing the same for every $i=1,\ldots,s$ we conclude that the image of $\pi_2:U(\omega_\LL)\rightarrow D(\omega_\LL)$ is $D(\omega_\LL,\overline{f})$.
\end{proof}

\section{The singular set}
\label{Isolated points of the singular set}
\hspace{9pt}In Section \ref{sec-5-1}, we use the decomposition of the singular set of a foliation $\omega\in\FF^1(\PP^n)(e)$, given by \cite{fmi}, and count the isolated points of $Sing(\omega)$ using the Hilbert polynomial of $\overline{\UU}(\omega)$, see Theorem \ref{teo-3}. In Section \ref{sec-5-2}, we show that the dimension of the classes of isomorphism projective unfoldings $\overline{U}(\omega)=\overline{\UU}(\omega)(e)$, does not succeed compute the number of isolated points of the singular set of $\omega$, by making some explicit computations.

\subsection{Counting isolated points of the singular set}\label{sec-5-1}
\hspace{9pt}Through this section we want to consider a foliation $\omega\in\FF^1(\PP^n)(e)$ of rational or logarithmic type. For that, we are going to extend the notation of a logarithmic form $\omega$ as
\begin{equation}\label{omega_log}
\omega = \sum_{i=1}^s \lambda_i\  F_i \ df_i 
\end{equation}
to the case where $s\geq 2$. We will say that $\omega$ is generic if $\omega$ is in the generic open sets $\mathcal{U}_\RR$ or $\mathcal{U}_\LL$.

Let us name $D_i$ the hypersurfaces defined by the functions $f_i$ and $D_{ij}$ the intersections $D_i\cap D_j$. We define the ideals $L_{ij}$ and $L = \bigcap_{i<j} L_{ij}$ associated to the varieties $D_{ij}$ and $Z = \bigcup D_{ij}$ respectively.

\

From \cite[Lemma 1.4, p. 8]{suwa-multiform} we have:

\prop\label{I=L} Let $\omega\in\mathcal{U}_\LL\subset\LL(n,\overline{d})$ be a generic logarithmic foliation. Then $I(\omega)=L$.

\

Let state the following result from \cite[Theorem, p.~3]{fmi}. Even if the authors focus in logarithmic foliations, there are no constrains to the case where $s=2$, which we consider here as well.
\prop\label{prop-fmi} Let $\omega\in \FF^1(\PP^n)(e)$ be a generic rational or logarithmic foliation. We can decompose the singular set of $\omega$ as the disjoint union
\[
Sing(\omega) = Z\cup Q
\]
where $Z = \bigcup_{i<j} D_{ij}$ and $Q$ is a set of finite points in $\PP^n$, consisting of $N(n,\overline{d})$ points counted with multiplicity. Even more, if any $d_i>1$ then $N(n,\overline{d})>0$.

\

Before stating our result, we need a technical definition.

\deff We are going to say that two graded $S$-modules $M$ and $N$ are stably isomorphic $M\simeq_s N$, if there exists $k_0\in\NN$ such that $M(k)\simeq N(k)$ for every $k\geq k_0$.

\teo\label{teo-3} Let $\omega\in\FF^1(\PP^n)(e)$ be a generic rational or logarithmic foliation. Then, the Hilbert polynomial $P_{\overline{\UU}(\omega)}$ of $\overline{\UU}(\omega)$ is constant and verifies
\[
P_{\overline{\UU}(\omega)} \equiv N(n,\overline{d}),
\]
where $N(n,\overline{d})$ is the number of isolated points of $Sing(\omega)$, counted with multiplicities.

\begin{proof}

Let us call $\OO_{Sing{(\omega)}}$, $\OO_Z$  and $\OO_Q$ to the structural sheafs of the correspondent varieties in $\PP^n$, from Proposition \ref{prop-fmi}. Being the union between $Z$ and $Q$ disjoint, we have the exact sequence of sheafs
\[
\xymatrix{
0 \ar[r] & \OO_Q \ar[r] & \OO_{Sing(\omega)} \ar[r] & \OO_Z \ar[r] & 0.
}
\]
Because of the annihilation of the higher cohomology of $\OO_Q$, applying the functor $\Gamma_*$ we get an exact sequence of graded $S$-modules
\begin{equation}\label{Gamma_*}
\xymatrix{
0 \ar[r] & \Gamma_*\OO_Q \ar[r] & \Gamma_*\OO_{Sing(\omega)} \ar[r] & \Gamma_*\OO_Z \ar[r] & 0.
}
\end{equation}

We can define another exact sequence of graded $S$-modules with the ideals $I(\omega)$ and $J(\omega)$ as 
\begin{equation}\label{suc-I/J}
\xymatrix{
0 \ar[r] & I(\omega)/J(\omega) \ar[r]^{} &  S/J(\omega) \ar[r]^{} & S/I(\omega) \ar[r] & 0.
}
\end{equation}

Writing $\omega$ as
\begin{equation*}\label{rac-log}
\omega = \sum_{i=1}^s\ \lambda_i\ F_i\ df_i = \sum_{i=0}^n \ A_i\ dx_i
\end{equation*}
we rapidly see that $J(\omega)=(A_0,\ldots, A_n)$, which implies $\Gamma_*\OO_{Sing(\omega)}\simeq_s S/J(\omega)$.
On the other side, by Propositions \ref{I=L} and \ref{prop-fmi} we have $I(\omega) = L = (F_1,\ldots,F_s)$ which implies $\Gamma_*\OO_Z\simeq_s S/I(\omega)$.

Comparing eq. \ref{Gamma_*} with eq. \ref{suc-I/J}, by the additivity of the Hilbert polynomial, we find the equalities
\[
\begin{aligned}
P_{\OO_Q} &= P_{Sing(\omega)} - P_{\OO_Z} = P_{S/J(\omega)} - P_{S/I(\omega)} = \\
&= P_{I(\omega)/J(\omega)}.
\end{aligned}
\]

Again, since $dim(\OO_Q)=0$ we have that $P_{\OO_Q}\equiv N(n,\overline{d})$ and by Proposition \ref{I/J} we get our result.
\end{proof}

\subsection{Some examples}\label{sec-5-2}

\hspace{9pt}We can reformulate Theorem \ref{teo-3} in the following way:

\teobis\label{cor5} Let $\omega\in\FF^1(\PP^n)(e)$ be a generic rational or logarithmic foliation. Then, there exits $a_\omega\in\NN$ such that, for $a\geq a_\omega$, we have
\[
dim_{\CC}\left(\overline{\UU}(\omega)(a)\right) = N(n,\overline{d}).
\]

We tried to see if the $a_\omega$ above, could be taken lower than the degree $e$ of $\omega$, to be able to compute $N(n,\overline{d})$ with $\overline{U}(\omega)$.
By explicit computations we found a negative answer for that.

For example, following the notation of eq. \ref{omega_log}, let us take $\omega_{(1,4)}\in\mathcal{U}_{\RR}\subset\RR(3,(1,4))$ defined with 
\[
f_1 = x_1\qquad \text{ and } \qquad f_2 = x_1^4+x_2^4+x_3^4+x_4^4
\]
and $\omega_{(1,2,3)}\in\mathcal{U}_{\LL}\subset\LL(3,(1,2,3))$ defined with 
\[
\begin{aligned}
&f_1 = x_1 \\
&f_2 = x_1^2-x_2^2+ix_3^2-ix_4^2\\
&f_3 = x_1^2x_2+x_1x_3^2+x_2^2x_4+x_3x_4^2 
\end{aligned}
\qquad
\begin{aligned}
 &\lambda_1=1 \ \ \ \lambda_2=i \\ &\lambda_3 = -\frac{2}{3}(1+i).
\end{aligned}
\]
With the help of a computer running \cite{diffalg} and \cite{M2}, we can compute the dimensions of $\overline{\UU}(\omega)(a)$ for enough $a\in\NN$ and see that they stabilize in $N(n,\overline{d})$, after the degree $e$ of the respective differential form.  We summarize that information in the following table:

\[
\begin{tabular}{c c c c c c c c c c }\label{tabla2}
& & & \multicolumn{4}{c}{$dim_{\CC}\left(\overline{\UU}(\omega)(a)\right)$} & & & \\[1.5ex]
\multicolumn{1}{c|}{$a$} & 1 & 2 & 3 & 4 & 5 & 6 & 7 & \ldots & \multicolumn{1}{||c}{$N(3,\overline{d})$ }\\
\multicolumn{1}{c|}{}& & & & & & & & & \multicolumn{1}{||c}{}\\[-2.0ex]
\cline{1-10}
\multicolumn{1}{c|}{}& & & & & & & & & \multicolumn{1}{||c}{}\\[-2.0ex]
\multicolumn{1}{c|}{$\omega_{(4,1)}$}   & 4 & 10 & 17 & 23 & \framebox{26} & 27 & 27 & \ldots & \multicolumn{1}{||c}{27}\\[0.2ex]
\multicolumn{1}{c|}{$\omega_{(1,2,3)}$} & 0 & 1  & 5  & 11 & 17 & \framebox{21} & 22 & \ldots  & \multicolumn{1}{||c}{22}\\
\end{tabular}
\]

\

\noindent where we write inside a box the dimension $dim_\CC\left(\overline{\UU}(\omega)(e)\right) = dim_\CC\left(\overline{U}(\omega)\right)$.

\section{Regularity}
\label{regularity}

\hspace{9pt}Along this section, let us fix $\omega\in\FF^1(\PP^n)(e)$ and regard it as an affine form in $\Omega^1_S$. We first extend the complex involved in the definition of regularity to a long complex $C^\bullet(\omega)$ of differential operators over $S$, see Definition \ref{complex_c}. Then, we prove that the cycles and borders of $C^\bullet(\omega)$ and $L^\bullet(\omega)$ are isomorphic, see Theorem \ref{triangle-dw}, relating the notion of regularity to the linear complex $L^\bullet(\omega)$. Finally, we show that the notion of regularity can be completely reinterpreted in terms of unfoldings, for the case of rational and logarithmic foliations, see Theorem \ref{teo-4}.

\subsection{The complex $C^\bullet({\omega})$}

\hspace{9pt}Let us recall the notion of regularity introduced in \cite[p.~17]{caln}:

\deff Let $\omega_0$ be an integrable, homogeneous differential 1-form in $\Omega^1_S(e)$. Then, $\omega_0$ it is said to be regular if for every $a<e$ the sequence
\begin{equation*}\label{complejo_cln}
\xymatrix@R=0pt{
T_S(a-e) \ar[r]^-{L(\omega_0)} & \Omega^1_{S}(a) \ar[r]^-{\omega_0\vartriangle} & \Omega^3_{S}(a+e)\\
X \ar@{|->}[r] & L_X(\omega_0) & &\\
&\eta \ar@{|->}[r] & \omega_0\vartriangle\eta & 
}
\end{equation*}
is exact in degree 1, {\it i.e.}, $Im(L_-(\omega_0)) = Ker(\omega_0\vartriangle- )$.

\

We are going to extend the definition of the differential operator $\omega_0\vartriangle-$ to the exterior algebra $\Omega_{S} = \bigoplus\limits_{r\geq 0} \Omega^r_{S}$. Let  $\tau\in\Omega^r_{S}$ and $\kappa(r):=\frac{r+1}{2}$, then we define 
\[
\omega_0\vartriangle\tau := \omega_0\wedge d\tau + \kappa(r) \ d\omega_0\wedge\tau.
\]

\prop \label{int-complejo} A differential form $\omega_0\in H^0(\Omega^1_{\PP^n}(e))$ is integrable if and only if $(\omega_0\vartriangle~)\circ(\omega_0\vartriangle~) \equiv 0$.
\begin{proof}
Let us take $\tau\in\Omega^r_{S}(p)$. Evaluating we get
\begin{multline}\label{ecu9}
\omega_0\vartriangle(\omega_0\vartriangle\tau)=\big(2(\kappa(r)+1)\big)\ \omega_0\wedge d\omega_0\wedge d\tau + \\ 
+ \ \big(\kappa(r)(\kappa(r)+1)\big) \ d\omega_0\wedge d\omega_0\wedge\tau.
\end{multline}
If we suppose that $\omega_0$ is integrable, from eq. \ref{ecu9} the first implication is clear.

For the other implication, take $\tau\in\Omega^r_{S}(p)$ and decompose it as $\tau = \tau_r + \tau_d$, using eq. \ref{lie-desc}. First equalize eq. \ref{ecu9} to 0, considering $\tau=\tau_r$ only.

By applying the exterior differential to $\omega_0\vartriangle(\omega_0\vartriangle\tau)$ and contracting with the radial field $R$, we see that
\[
2e \ \omega_0\wedge d\omega_0\wedge d\tau_r + p\ d\omega_0\wedge d\omega_0\wedge\tau_r = 0 .
\]
If we choose $p$ such that the $2$-uplas of coefficients $(2(\kappa(r)+1),\kappa(r)(\kappa(r)+1))$ and $(2e,p)$ are linearly independent, we can cancel terms and see
\[
d\omega_0\wedge d\omega_0 \wedge \tau_r =0.
\]
Now equalize eq. \ref{ecu9} to 0 considering $\tau=\tau_d$. We immediately get $d\omega_0\wedge d\omega_0\wedge \tau_d = 0$.
This way we see that $d\omega_0\wedge d\omega_0\wedge \tau = 0$ for every $\tau\in\Omega^r_S(p)$ and the result follows.
\end{proof}

\

For $\omega\in\FF^1(\PP^n)(e)$, the above property allows us to define the complex of $\CC$-vector spaces $C^\bullet(\omega)$ in the following way.

\deff\label{complex_c} We define the graded complex  $C^\bullet(\omega)$ associated to $\omega$, as
\[
\xymatrix{
C^\bullet(\omega):& T_S\ar[r]^-{\omega\vartriangle} & \Omega^1_{S} \ar[r]^-{\omega\vartriangle} & \Omega^3_{S} \ar[r]^-{\omega\vartriangle} & \ldots & &
}
\]
where $C^s=\Omega^{2s-1}_S$ for $s\geq 0$ and the 0-th differential is defined as $\omega\vartriangle X:=L_X(\omega)= i_Xd\omega + di_X\omega$. The grading of $C^\bullet(\omega)$ is given by the decomposition $C^\bullet(\omega) = \bigoplus_{a\in\NN} C^\bullet(\omega)(a)$, where $C^\bullet(\omega)(a)$ is the complex of finite vector spaces
\[
\xymatrix@C=35pt{
C^\bullet(\omega,a):& T_S(a-e)\ar[r]^-{\omega\vartriangle} & \Omega^1_{S}(a) \ar[r]^-{\omega\vartriangle} & \Omega^3_{S}(a+e) \ar[r]^-{\omega\vartriangle} & \ldots .
}
\]

\

\remo Let $d:M\to M$ be a $\CC$-linear function, and $M$ an $S$-module. We say that $d$ is a differential operator (of order 1) over $S$ if, for every $f\in S$, the application $m\mapsto d(f.m)-f.d(m)$ is $S$-linear on $M$. It is immediate to check that the complex $C^\bullet(\omega)$ is a complex of differential operators over $S$.

\

Since $C^\bullet(\omega)$ is not $S$-linear, we will not find a morphism of complexes between $C^\bullet(\omega)$ and $L^\bullet(\omega)$. Anyway, we will be able to find $\CC$-linear isomorphisms on every degree, but one. For being able to compare these two complexes we need the following technical elements:

\deff \label{formula} We are going to say that a 5-uple of indexes $(r,s,p,a,e)\in\ZZ^5$ is {\it admissible}, if $a,p\in\NN$, $e\geq 2$, $r\geq -1$ and
\begin{equation*}\label{ec-formula}
\left\{
\begin{aligned}
s & = \kappa(r)\\
p &= e(s-1)+ a = e\left(\frac{r-1}{2}\right) + a.
\end{aligned}
\right.
\end{equation*}

\lemma The $5$-uple of indexes $(r,s,p,a,e)\in\ZZ^5$ is admissible  if and only if the following equalities hold
\[
C^s(\omega,a)=L^s(\omega,a)=\Omega^{r}_{S}(p).
\]

\deff Let us consider an admissible uple of indexes $(r,s,p,a,e)\in\ZZ^5$. We define the family of graded linear maps $\left\{\varphi^s_a:\Omega^r_{S}\to\Omega^r_{S}\right\}$, such that in each homogeneous component of degree $p$, $\varphi^s_a$ is defined by
\[
\xymatrix@R=5pt@C=70pt{
\Omega^r_{S}(p) \ar[r]^-{\varphi^s_a} & \Omega^{r}_{S}(p) \\
\save[] *{\tau = \tau_r + \tau_d\ }="s" &  \save[] +<15pt,0pt> *{\ \big(e\kappa(r)-p\big)\tau_r + e\kappa(r)\tau_d}="t" \restore \\ \restore \ar@{|->}"s";"t"
}
\]
for $s\geq 1$. For $s=0$, $\varphi^s_a: T_S \to T_S $ is the identity map.

\

\teo\label{triangle-dw} Let  $\omega\in\FF^1(\PP^n)(e)$ and $(r,s,p,a,e)\in\ZZ^5$ admissible, such that $a\neq e$. Then, the family $\left\{\varphi^s_a\right\}$ induces isomorphisms of $\CC$-vector spaces
\[
\mathcal{Z}^s(C^\bullet(\omega))\simeq \mathcal{Z}^s(L^\bullet(\omega))\qquad \text{and}\qquad \mathcal{B}^s(C^\bullet(\omega))\simeq \mathcal{B}^s(L^\bullet(\omega)).
\]
\begin{proof}
Since $a\neq e$, $e\kappa(r)-p\neq 0$ and then the $\varphi^s_a$ are all isomorphisms over $\Omega^r_S(p)$. We just need to show that $\varphi^s_a$ sends isomorphically the kernel of $\omega\vartriangle-$ to the kernel of $d\omega\wedge-$.

Let us take $\tau\in\Omega^r_{S}(p)$ such that $\omega\vartriangle\tau=0$. Decomposing $\tau=\tau_r+\tau_d$ as in eq. \ref{lie-desc}, we have
\begin{equation}\label{equ-1}
\omega\vartriangle\tau = \omega\wedge d\tau_r + \kappa(r)\ d\omega\wedge\tau_r + \kappa(r)\ d\omega\wedge\tau_d = 0.
\end{equation}
Applying exterior differential and contracting with $R$ we get
\begin{equation}\label{equ-2}
e\ \omega\wedge d\tau_r + p\ d\omega\wedge\tau_r = 0.
\end{equation}
Operating with eqs. \ref{equ-1} and \ref{equ-2}, we can write 
\[
d\omega\wedge \Big(\big(e\kappa(r)-p\big) \tau_r + e\kappa(r) \tau_d\Big) = d\omega\wedge\varphi^s_a(\tau) = 0
\]
from where we get one implication.

\

Let us suppose now that $d\omega\wedge\tau = 0$. Differentiating and contracting as before, we end up noticing that eq. \ref{equ-2} still holds. If we apply the inverse function to $\tau$
\[
(\varphi^s_a)^{-1}(\tau) = \frac{1}{e\kappa(r)-p}\tau_r + \frac{1}{e\kappa(r)}\tau_d
\]
and compose with $\omega\vartriangle-$, we find the expression
\begin{equation}\label{equ-3}
\omega\wedge\left(\frac{1}{e\kappa(r)-p}d\tau_r\right)+\frac{\kappa(r)}{e\kappa(r)-p}d\omega\wedge\tau_r + \frac{1}{e}d\omega\wedge\tau_d.
\end{equation}
Since $d\omega\wedge d\tau_r = - d\omega\wedge\tau_d$, we can simplify eq. \ref{equ-3} and see that
\[
e\ \omega\wedge d\tau_r + e\kappa(r)\  d\omega\wedge\tau_r - \big(e\kappa(r)-p\big) d\omega\wedge\tau_r
\]
equals the right side of eq. \ref{equ-2}. This way, we conclude that  $\omega\vartriangle \Big(\left(\varphi^s_a\right)^{-1}(\tau)\Big) = 0$.

\

For the case where $r=-1$, let us consider $X\in T_S(b)$, with $b=a-e\neq 0$. Suppose
\begin{equation}\label{equ-4}
\omega\vartriangle X= i_Xd\omega + di_X\omega=0.
\end{equation}
Contracting with $R$ we get $i_X\omega=0$, which implies $di_X\omega = 0$. Together with eq. \ref{equ-4} we see that $d\omega\wedge X=i_Xd\omega = 0$.

Finally, suppose that $i_Xd\omega = 0$. Contracting with $R$ we get the other necessary term to obtain the formula $\omega\vartriangle X=0$.
\end{proof}

\cor\label{cor-C-L} For every $a\neq e$ we have
\[
dim_{\CC}\left(H^s(C^\bullet(\omega,a)\right) \ =\ dim_{\CC}\left(H^s(L^\bullet(\omega,a))\right) .
\]

\

Applying Corollary \ref{cor-C-L} above, we can state our final result and relate regularity to first order unfoldings.

\teo\label{teo-4} Let $\omega\in\FF^1(\PP^n)(e)$ be a generic rational or logarithmic foliation. Then, $\omega\in\Omega^1_S(e)$ is regular if and only if $\overline{U}(\omega)=0$.
\begin{proof}
Using Corollary \ref{H^1I/J} and Corollary \ref{cor-C-L} we have that $\omega$ is regular if and only if $\left(I(\omega)/J(\omega)\right)(a)=0$ for every $a<e$.

If $\omega$ is of type $(1,\ldots,1)$ we rapidly see that $I(\omega)=J(\omega)$.

By the definition of $J(\omega)$ and by Theorem \ref{I-meromorphic} or Theorem \ref{I-multiform}, we now that $I(\omega)$ and $J(\omega)$ are generated in degrees lower than $e$. Suppose $\omega$ is of type $(d_1,\ldots,d_s)$ and some $d_k>1$. Then, $\omega$ is not regular since $I(\omega)(e-d_k)\neq 0$ and $J(\omega)(e-d_k)=0$.

Putting together Corollary \ref{H^1I/J}, Proposition \ref{prop-fmi} and Theorem \ref{teo-3},  we see that $N(n,\overline{d})=P_{I(\omega)/J(\omega)}\neq 0$ what forces $(I(\omega)/J(\omega))(e)=\overline{U}(\omega)$ to be $\neq 0$.
\end{proof}

By the above proof we can also claim:

\cor Let $\omega\in\FF^1(\PP^n)(e)$ be a generic rational or logarithmic foliation of type $\overline{d}=(d_1,\ldots,d_s)$ for $s\leq n+1$. Then, $\omega\in\Omega^1_S(e)$ is regular if and only if $\overline{d}=(1,\ldots,1)$.

\begin{appendices}

\section{}

\hspace{9pt}\label{proof2}
We can make an analytic proof of Theorem \ref{teo-1} without the need of the classification of $D(\omega_\RR)$ of Theorem \ref{tangent-rational}:

\begin{proof}[Alternative proof of Theorem \ref{teo-1}]
Let us consider $\eta\in H^0(\Omega^1_{\PP^n}(e))$ such that verifies the equation 
\begin{equation}\label{equ}
\omega_{\RR}\wedge d\eta + d\omega_{\RR}\wedge \eta = 0.
\end{equation}

Let us pullback $\omega_\RR$ and $\eta$ to the affine space $\CC^{n+1}$ and take a point $p$ such that 
\[
 d\omega_{\RR}(p)=(r+s) df_1(p)\wedge d{f_2}(p)\neq 0.
\]
Then functions $f_1, {f_2}$ are transversal in a neighborhood $W$ of $p$ and we can choose a coordinate system of the form $(W,\overline{\varphi})$ such that $\overline{\varphi} = (f_1,{f_2},\varphi_1,\ldots,\varphi_{\ell})$.

\

In this neighborhood, $\eta$ can be written as
\[
\eta = h_{f_1} \ d{f_1} + h_{f_2} \ d{f_2} + \sum_{i=1}^\ell h_i \ d\varphi_i
\]
and $d\eta$ as
\begin{align*}
d\eta &= \left( \frac{\partial h_{f_2}}{\partial {f_1}} - \frac{\partial h_{f_1}}{\partial {f_2}} \right) d{f_1}\wedge d{f_2} + \sum_{i=1}^\ell \left( \frac{\partial h_i}{\partial {f_1}} - \frac{\partial h_{f_1}}{\partial \varphi_i}\right) d{f_1}\wedge d\varphi_i\  + \\
&\hspace{1cm} +\sum_{i=1}^\ell \left( \frac{\partial h_i}{\partial {f_2}} - \frac{\partial h_{f_2}}{\partial \varphi_i}\right) d{f_2}\wedge d\varphi_i + \sum_{\substack{i,j=1\\i<j}}^\ell \left( \frac{\partial h_j}{\partial \varphi_i} - \frac{\partial h_i}{\partial \varphi_j}\right) d\varphi_i\wedge d\varphi_j.
\end{align*}

Using the formulas above we can expand eq. \ref{equ} and get
\begin{align*}
&\omega_{\RR}\wedge d\eta  + d\omega_{\RR}\wedge\eta = \\
&=\ \sum_{i=1}^\ell \left[ r {f_1} \ \left( \frac{\partial h_i}{\partial {f_1}} - \frac{\partial h_{f_1}}{\partial \varphi_i}\right) + s {f_2}\ \left( \frac{\partial h_i}{\partial {f_2}} - \frac{\partial h_{f_2}}{\partial \varphi_i}\right) \right]d{f_1}\wedge d{f_2}\wedge d\varphi_i \ + \\
% \end{align*}
% \begin{align*}
&\hspace{2cm}+  \sum_{\substack{i,j=1\\i<j}}^\ell r{f_1}\left( \frac{\partial h_j}{\partial \varphi_i} - \frac{\partial h_i}{\partial \varphi_j}\right) d{f_2}\wedge d\varphi_i \wedge d\varphi_j \ + \\
% \end{align*}
% \begin{align*}
&\hspace{2cm}-  \sum_{\substack{i,j=1\\i<j}}^\ell s{f_2}\left( \frac{\partial h_j}{\partial \varphi_i} - \frac{\partial h_i}{\partial \varphi_j}\right) d{f_1}\wedge d\varphi_i\wedge d\varphi_j = 0.
\end{align*}

From the last two summations we obtain the equality $\frac{\partial h_j}{\partial \varphi_i} =\frac{\partial h_i}{\partial \varphi_j}$ for every $i,j$.

Now, let us take $k\in\{1,\ldots,\ell\}$ and define the primitive function  $h = \int h_kd\varphi_k$. We can compute a partial derivative of $h$ in $W$ as
\[ 
\frac{\partial h}{\partial \varphi_i} = \int \frac{\partial h_k}{\partial \varphi_i} \ d\varphi_k = \int \frac{\partial h_i}{\partial \varphi_k} \ d\varphi_k = h_i
\]
and then, express the differential of $h$ as
\[
dh = \frac{\partial h}{\partial {f_1}} \ d{f_1} + \frac{\partial h}{\partial {f_2}} \ d{f_2}+ \sum_{i=1}^\ell h_i \ d\varphi_i.
\]

If we consider the 1-form
\[
\eta - dh = \left( h_{f_1} - \frac{\partial h}{\partial {f_1}} \right) d{f_1} + \left( h_{f_2} - \frac{\partial h}{\partial {f_2}} \right) d{f_2}
\]
we clearly have $d\omega_{\RR}\wedge(\eta-dh)=0$.

Both $\omega_\RR$ and $\eta$ are homogeneous of degree $e$ and descend to projective space. Let us restrict to the homogeneous part of degree $e$ of the previous equation and call $h_e$ the homogeneous component of $h$ of that degree. We can contract the equation $d\omega_\RR\wedge(\eta-dh_e)=0$ with the radial field $R$ and get
\[
 h_ed\omega = \omega\wedge(\eta-dh_e)
\]
showing that the pair $(h_e,\eta)\in U(\omega_\RR)$ and projects to the deformation $\eta$.
\end{proof}

\end{appendices}

% % % % % % % % % % % % 
% % % % % % % % % % % % 
% BIBLIOGRAPHY
% % % % % % % % % % % %
% % % % % % % % % % % %

\bibliography{biblio}{}

\begin{thebibliography}{GMLN91}

\bibitem[CA94]{omegar}
O.~Calvo-Andrade.
\newblock Irreducible components of the space of holomorphic foliations.
\newblock {\em Math. Ann.}, 299(4):751--767, 1994.

\bibitem[CGM]{fjc}
F.~Cukierman, J.~Gargiulo, and C.~Massri.
\newblock On the stability of logarithmic differential one-forms.
\newblock To appear.

\bibitem[CLN82]{caln}
C.~Camacho and A.~Lins~Neto.
\newblock The topology of integrable differential forms near a singularity.
\newblock {\em Inst. Hautes \'Etudes Sci. Publ. Math.}, (55):5--35, 1982.

\bibitem[CLN96]{celn}
D.~Cerveau and A.~Lins~Neto.
\newblock Irreducible components of the space of holomorphic foliations of
  degree two in {$\mathbf C{\rm P}(n)$}, {$n\geq 3$}.
\newblock {\em Ann. of Math. (2)}, 143(3):577--612, 1996.

\bibitem[CPV09]{fji}
F.~Cukierman, J.~V. Pereira, and I.~Vainsencher.
\newblock Stability of foliations induced by rational maps.
\newblock {\em Ann. Fac. Sci. Toulouse Math. (6)}, 18(4):685--715, 2009.

\bibitem[CSV06]{fmi}
F.~Cukierman, M.~G. Soares, and I.~Vainsencher.
\newblock Singularities of logarithmic foliations.
\newblock {\em Compos. Math.}, 142(1):131--142, 2006.

\bibitem[Die81]{dieudonne}
J.~Dieudonn{\'e}.
\newblock {\em \'{E}l\'ements d'analyse. {T}ome {I}}.
\newblock Cahiers Scientifiques [Scientific Reports], XXVIII. Gauthier-Villars,
  Paris, third edition, 1981.
\newblock Fondements de l'analyse moderne. [Foundations of modern analysis],
  Translated from the English by D. Huet, With a foreword by Gaston Julia.

\bibitem[DMM]{diffalg}
M.~Dubinsky, C.~D. Massri, and A.~Molinuevo.
\newblock diff{A}lg, a differential algebra library.
\newblock \\ Available at \href{https://savannah.nongnu.org/projects/diffalg/}%
  {https://savannah.nongnu.org/projects/diffalg/}.

\bibitem[Eis95]{eisenbud}
D.~Eisenbud.
\newblock {\em Commutative algebra}, volume 150 of {\em Graduate Texts in
  Mathematics}.
\newblock Springer-Verlag, New York, 1995.
\newblock With a view toward algebraic geometry.

\bibitem[GMLN91]{gmln}
X.~G{\'o}mez-Mont and A.~Lins~Neto.
\newblock Structural stability of singular holomorphic foliations having a
  meromorphic first integral.
\newblock {\em Topology}, 30(3):315--334, 1991.

\bibitem[GS]{M2}
D.~R. Grayson and M.~E. Stillman.
\newblock Macaulay2, a \mbox{software} system for research in algebraic
  geometry.
\newblock \\ Available at \href{http://www.math.uiuc.edu/Macaulay2/}%
  {http://www.math.uiuc.edu/Macaulay2/}.

\bibitem[Har77]{hart}
R.~Hartshorne.
\newblock {\em Algebraic geometry}.
\newblock Springer-Verlag, New York, 1977.
\newblock Graduate Texts in Mathematics, No. 52.

\bibitem[Jou79]{jou}
J.~P. Jouanolou.
\newblock {\em \'{E}quations de {P}faff alg\'ebriques}, volume 708 of {\em
  Lecture Notes in Mathematics}.
\newblock Springer, Berlin, 1979.

\bibitem[Mal76]{m-f1}
B.~Malgrange.
\newblock Frobenius avec singularit\'es. {I}. {C}odimension un.
\newblock {\em Inst. Hautes \'Etudes Sci. Publ. Math.}, (46):163--173, 1976.

\bibitem[Mal77]{m-f2}
B.~Malgrange.
\newblock Frobenius avec singularit\'es. {II}. {L}e cas g\'en\'eral.
\newblock {\em Invent. Math.}, 39(1):67--89, 1977.

\bibitem[Mum99]{mumford}
D.~Mumford.
\newblock {\em The red book of varieties and schemes}, volume 1358 of {\em
  Lecture Notes in Mathematics}.
\newblock Springer-Verlag, Berlin, expanded edition, 1999.
\newblock Includes the Michigan lectures (1974) on curves and their Jacobians,
  With contributions by Enrico Arbarello.

\bibitem[Suw83a]{suwa-multiform}
T.~Suwa.
\newblock Unfoldings of foliations with multiform first integrals.
\newblock {\em Ann. Inst. Fourier (Grenoble)}, 33(3):99--112, 1983.

\bibitem[Suw83b]{suwa-meromorphic}
T.~Suwa.
\newblock Unfoldings of meromorphic functions.
\newblock {\em Math. Ann.}, 262(2):215--224, 1983.

\bibitem[Suw95]{suwa-review}
T.~Suwa.
\newblock Unfoldings of codimension one complex analytic foliation
  singularities.
\newblock In {\em Singularity theory ({T}rieste, 1991)}, pages 817--865. World
  Sci. Publ., River Edge, NJ, 1995.

\bibitem[War83]{warner}
F.~W. Warner.
\newblock {\em Foundations of differentiable manifolds and {L}ie groups},
  volume~94 of {\em Graduate Texts in Mathematics}.
\newblock Springer-Verlag, New York, 1983.
\newblock Corrected reprint of the 1971 edition.

\end{thebibliography}
% \addcontentsline{toc}{section}{References}
\bibliographystyle{alpha}

\

\address{\sc Ariel Molinuevo {\hfill{\sf amoli@dm.uba.ar}} \\Departamento de Matem\'atica, FCEyN\\Universidad de Buenos Aires\\Ciudad Universitaria, Pabell\'on I\\ CP C1428EGA\\ Buenos Aires\\ Argentina}

\end{document}